\documentclass[10pt]{article}
\usepackage{amsmath,amsbsy,amsfonts,amssymb}

\usepackage[frenchb]{babel}

\begin{document}

\title{Classification  des s\'eries discr\`etes pour certains groupes classiques $p$-adiques}
\author{C. M{\oe}glin}
\date{}
\maketitle
\rightline{En l'honneur de Roger Howe pour son 60e anniversaire}

\

L'article a pour but de montrer que seule l'hypoth\`ese que des lemmes fondamentaux pour des paires bien 
pr\'ecises de groupes est suffisante pour avoir les hypoth\`eses de \cite{europe} et \cite{ams} permettant de classifier les s\'eries discr\`etes des groupes symplectiques et orthogonaux ''impairs'' p-adiques.
Le point de d\'epart est l'id\'ee d'Arthur exprim\'ee pleinement sous forme de conjecture en \cite{arthur} et 
sous forme de r\'esultats modulo des lemmes fondamentaux en \cite{arthurnouveau}.

Le corps de base, not\'e $F$, est un corps $p$-adique; on suppose que $p\neq 2$ car nos r\'ef\'erences, surtout \cite{waldspurger1} utilisent cette hypoth\`ese. Il n'intervient pas dans notre propre travail.

On consid\`ere les points sur $F$ d'un groupe  classique, c'est-\`a-dire le groupe des automorphismes d'une forme orthogonale ou symplectique sur un $F$-espace vectoriel, $V$, de dimension finie; ici on se limite aux espaces $V$ symplectiques ou orthogonaux de dimension impaire.  On note $n$ est le rang du groupe, c'est-\`a-dire la partie enti\`ere, $[dim\, V/2]$ et $G(n)$ ce groupe p-adique. Si $V$ est orthogonal, on ne suppose pas $V$ maximalement 
d\'eploy\'ee  pour obtenir le r\'esultat principal toutefois le cas non maximalement d\'eploy\'e s'obtient \`a la fin \`a partir du cas maximalement d\'eploy\'e et dans cette introduction on suppose $G(n)$ 
d\'eploy\'e. On note $G^*(n)$ le groupe complexe des automorphismes de la forme duale au sens de Langlands,  c'est-\`a-dire si $V$ est orthogonal de dimension impaire $G^*=Sp(2n,{\mathbb C})$, si $V$ est symplectique de dimension paire $G^*=SO(2n+1,{\mathbb C})$.
On notera $n^*$ l'entier dimension de la repr\'esentation naturelle de $G^*$; c'est-\`a-dire que $n^*=2n$ si $G^*(n)$ est un groupe symplectique et $2n+1$ si $G^*(n)$ est un groupe orthogonal. Quand un groupe lin\'eaire $GL(m,F)$ ($m\in {\mathbb N}$) est fix\'e, on note $\theta$ l'automorphisme de ce groupe 
$g\mapsto J_{m}^{-1}\, ^tg^{-1} J_{m}$, o\`u $J_{m}$ est essentiellement la matrice antidiagonale (cf \ref{notations}) et $\tilde{GL}(m,F)$ le produit semi-direct de $GL(m,F)$ par le groupe $\{1,\theta\}$.

L'id\'ee de base d'Arthur pour \'etudier les groupes classiques est que $G(n)$ est un groupe endoscopique principal pour l'ensemble $\tilde{G}(n^*)$ 
qui est par d\'efinition la composante connexe de $\theta$ dans le produit semi-direct $\tilde{GL}(n^*,F)$. C'est le groupe qui permet de faire un transfert entre les objects stables pour $G(n)$ et les objets stables pour $\tilde{G}(n^*)$. Arthur a aussi expliqu\'e comment obtenir un transfert pour les groupes orthogonaux sur un espace de dimension paire; la ''seule'' diff\'erence est que le transfert fait intervenir des facteurs de transfert que l'auteur ne ma\^{\i}trise pas.

L'hypoth\`ese faite dans ce papier est que le lemme fondamental pour ce couple soit vrai pour tout $m\leq n$, c'est indispensable pour l'\'enonc\'e m\^eme des r\'esultats. Mais il faut en plus le lemme fondamental pour tous les groupes endoscopiques  intervenant pour la stabilisation de $\tilde{G}(m^*)$, pour tout $m\leq n$. Et il y a un r\'esultat d'Arthur qui doit m\^eme, je crois, supposer la validit\'e des lemmes fondamentaux comme ci-dessus pour $m>n$; on peut sans doute borner $m$ mais sans doute pas suffisamment pour que des cas particuliers puissent vraiment \^etre d\'emontr\'es sans que le cas g\'en\'eral le soit. On n'a par contre pas besoin de leur analogue pond\'er\'e (sauf erreur); on utilise de fa\c{c}on d\'eterminante les r\'esultats locaux de  \cite{arthurnouveau} mais  pas les r\'esultats globaux.  

Pour la d\'etermination des repr\'esentations cuspidales et donc la classification compl\`ete, on a aussi besoin des lemmes fondamentaux pour tous les groupes endoscopiques elliptiques de $\tilde{GL}(n^*,F)$. Pour cette classification compl\`ete, nous nous sommes limit\'es aux groupes $G$ orthogonaux sur un espace de dimension impaire bien que les m\'ethodes soient g\'en\'erales; la raison de cette limitation est que l'on ne veut pas ici discuter le cas des groupes orthogonaux pairs qui m\'eritent une discussion \`a eux tout seuls. Pour eux, dans tous mes travaux, \`a la suite d'Adams, j'ai choisi de prendre comme groupe dual un groupe non connexe, le groupe orthgonal complexe; ce n'est pas le choix d'Arthur et il y a donc \`a \'ecrire les \'equivalences entre les 2 points de vue.

Une repr\'esentation dont la classe d'isomorphie est invariante sous l'action de $\theta$ sera dite une repr\'esentation $\theta$-invariante.
Pour tout $m$,  on appelle s\'erie $\theta$-dicr\`ete une repr\'esentation irr\'eductible temp\'er\'ee de $GL(m,F)$,  $\theta$-invariante  et qui n'est induite propre d'aucune repr\'esentation $\theta$-invariante d'un sous-groupe de Levi $\theta$-invariant de $GL(m,F)$.

Cette d\'efinition est beaucoup plus simple du c\^ot\'e groupe dual; on rappelle que la conjecture de Langlands locale pour les $GL$ (\cite{harris}, \cite{henniart}) et les  r\'esultats de Zelevinsky (\cite{zelevinsky}) permettent d'associer \`a toute 
repr\'esentation irr\'eductible temp\'er\'ee de $GL(m,F)$ une classe de conjugaison d'un morphisme continu, born\'e de $W_{F}\times SL(2,{\mathbb C})$ dans $GL(m,{\mathbb C})$. Pour $\pi^{GL}$ une repr\'esentation irr\'eductible temp\'er\'ee de $GL(m,F)$ on note $\psi_{\pi^{GL}}$, ce morphisme; le fait que $\pi^{GL}$ est $\theta$-invariante est exactement \'equivalent \`a ce que (quitte \`a conjuguer) $\theta \circ \psi_{\pi^{GL}}\circ \theta=\psi_{\pi^{GL}}$, o\`u  ici $\theta$ est l'automorphisme $g\mapsto \, ^tg^{-1}$ de $GL(m,{\mathbb C})$. Le fait que $\pi^{GL}$ est $\theta$-discr\`ete est \'equivalent \`a ce que le sous-groupe du centralisateur de $\psi_{\pi^{GL}}$ dans $GL(m,{\mathbb C})$ form\'e par les \'el\'ements $\theta$-invariants est fini.

Soit encore $m$ un entier et soit $\pi^{GL}$ une repr\'esentation temp\'er\'ee de $GL(m^*,F)$, $\theta$-invariante; on fixe un prolongement de $\pi^{GL}$ \`a $\tilde{GL}(m^*,F)$. Le caract\`ere de cette 
repr\'esentation restreint \`a $\tilde{G}(n^*)$ n'est pas n\'ecessairement stable. C'est m\^eme un r\'esultat profond de \cite{arthurnouveau} 30.1 de pr\'eciser que ce caract\`ere  est stable pr\'ecis\'ement quand le morphisme associ\'e se factorise par $G^*(n)$; c'est pour pouvoir utiliser ce r\'esultat que l'on a 
suppos\'e la totalit\'e des lemmes fondamentaux.  Et en suivant Arthur, quand le caract\`ere de $\pi^{GL}$ vue comme distribution sur $\tilde{G}(n^*)$ est stable, on d\'efinit le paquet de $\pi^{GL}$, not\'e $\Pi(\pi^{GL})$ comme \'etant l'ensemble  des repr\'esentations irr\'eductibles, $\pi$ de $G(m)$ tel que pour un bon choix de coefficients $c_{\pi}\in {\mathbb R}_{\geq 0}$ la distribution $tr \pi^{GL}$ de $GL(m^*,F)\times \theta$ soit un transfert de $\sum_{\pi\in \Pi(\pi^{GL})}c_{\pi}tr \pi$ (cf \ref{appartenance});  cette condition d\'efinit uniquement les 
\'el\'ements de $\Pi(\pi^{GL})$.   Si $\psi$ est le morphisme correspondant \`a $\pi^{GL}$, on note $\Pi(\psi)$ au lieu de $\Pi(\pi^{GL})$. Comme Arthur a en vue des r\'esultats globaux, il prolonge cette d\'efinition pour toute repr\'esentation $\pi^{GL}$ composante locale d'une forme automorphe cuspidale de $GL(m^*)$; la g\'en\'eralisation est \'evidente, c'est une induction, il prolonge d'ailleurs aussi les constructions pour obtenir toutes les formes automorphes de carr\'e int\'egrable mais nous n'en avons pas besoin ici.  Une remarque  par exemple dans  \cite{waldspurger1} VI.1 (ii)  est que si $\pi^{GL}$ est $\theta$-discr\`ete alors $\Pi(\pi^{GL})$ est form\'e de repr\'esentations elliptiques; comme les coefficients $c_{\pi}$ sont positifs d'apr\`es Arthur, on peut encore remplacer elliptiques par s\'eries discr\`etes. Toutefois, on red\'emontrera ce dont on a besoin pour \'eviter cet argument crois\'e.

Pour pouvoir utiliser ce qui pr\'ec\`ede,
la remarque \`a v\'erifier est que si $\pi$ est une s\'erie discr\`ete de $G(n)$ irr\'eductible alors le 
caract\`ere distribution de cette repr\'esentation appartient \`a un paquet stable, c'est-\`a-dire qu'il existe $\pi^{GL}$, $\theta$-discret comme ci-dessus tel que $\pi\in \Pi(\pi^{GL})$; cela est fait en \ref{appartenance}. Ensuite, il n'est pas difficile de v\'erifier que cela entra\^{\i}ne que le support cuspidal de $\pi$ est ''demi-entier'' (cf. \ref{definitionsupportcuspidal}) et par voie de cons\'equence, cela n\'ecessite que les points de r\'eductibilit\'e des induites de cuspidales soient demi-entiers (cf. \ref{pointsdereductibilite}). Ceci est une partie des 
hypoth\`eses de \cite{europe} et \cite{ams} et on peut aller plus loin pour avoir toutes les hypoth\`eses. Ainsi on associe \`a $\pi$ le morphisme $\psi_{\pi^{GL}}$ tel que $\pi\in \Pi(\psi_{\pi^{GL}})$. A ce point, il me semble que l'unicit\'e de $\pi^{GL}$ n'est pas compl\`etement claire;   on (re)trouvera cette unicit\'e de fa\c{c}on d\'etourn\'ee \`a la fin du travail en \ref{compatibilite}.

On peut  faire cette construction  pour $\pi_{0}$ une repr\'esentation cuspidale de $G(n)$ et donc obtenir un morphisme $\psi_{0}$ de $W_{F}\times SL(2,{\mathbb C})$ dans $GL(n^*,{\mathbb C})$.  Ce morphisme d\'etermine et est d\'etermin\'e par les propri\'et\'es d'irr\'eductibilit\'e des induites de la forme $St(\rho,a)\times \pi_{0}$, o\`u $\rho$ est une repr\'esentation irr\'eductible, cuspidale autoduale d'un groupe $GL(d_{\rho},F)$ et $a$ est un entier; $St(\rho,a)$ est la repr\'esentation de Steinberg g\'en\'eralis\'ee de $GL(ad_{\rho},F)$, l'induite \'etant une repr\'esentation de $G(n+ad_{\rho})$. Ce lien se fait exactement comme conjectur\'e en \cite{europe} et \cite{ams}:

on voit $\psi_{0}$ comme une repr\'esentation de $W_{F}\times SL(2,{\mathbb C})$ dans 
$GL(n_{0}^*,{\mathbb C})$. Toute repr\'esentation irr\'eductible de $W_{F}\times SL(2,{\mathbb C})$ est le produit tensoriel d'une repr\'esentation irr\'eductible, not\'ee $\rho$, de $W_{F}$, et d'une repr\'esentation de dimension finie $\sigma_{a}$ de $SL(2,{\mathbb C})$ uniquement d\'etermin\'ee par sa dimension not\'ee ici $a\in {\mathbb N}$; on note $d_{\rho}$ la dimension de la repr\'esentation $\rho$ et on identifie $\rho$ \`a une repr\'esentation cuspidale irr\'eductible de $GL(d_{\rho},F)$ gr\^ace \`a la correspondance locale de Langlands d\'emontr\'ee en \cite{harris}, \cite{henniart} et on garde la m\^eme notation $\rho$.  Et on montre qu'une telle repr\'esentation irr\'eductible $\rho\otimes \sigma_{a}$ de $W_{F}\times SL(2,{\mathbb C})$ intervient dans $\psi_{0}$ si et seulement si les 2 conditions suivantes sont r\'ealis\'ees:

\

(1)  $\rho\otimes \sigma_{a}$ \`a conjugaison pr\`es est \`a valeurs dans un groupe de m\^eme  type  que $G^*$, c'est-\`a-dire orthogonal si ce groupe est orthogonal et symplectique sinon;

(2)  l'induite pour le groupe $G(n_{0}+d_{\rho})$, 
$
St(\rho,a)\times \pi_{0}$ est irr\'eductible, o\`u $St(\rho,a)$ est la repr\'esentation de Steinberg bas\'ee sur $\rho$ du groupe $GL(ad_{\rho},F)$.

\

\noindent
Et on montre aussi que si $\rho\otimes \sigma_{a}$ intervient comme sous-repr\'esentation de $\psi_{0}$, elle y intervient avec multiplicit\'e 1. 

Ces propri\'et\'es entra\^{\i}nent l'unicit\'e de $\psi_{0}$ puisque l'on conna\^{\i}t sa d\'ecomposition en repr\'esentations irr\'educ\-tibles gr\^ace \`a des propri\'et\'es de la repr\'esentation $\pi_{0}$; c'est ce que l'on a appel\'e la connaissance des blocs de Jordan de $\pi_{0}$ et de $\psi_{0}$.

Le lien avec les points de r\'eductibilit\'e des induites de la forme $\rho\vert\,\vert^x\times \pi_{0}$, o\`u $\rho$ est comme ci-dessus et $x$ est un nombre r\'eel,  est le suivant. Pour $\rho$ fix\'ee vue comme une repr\'esentation irr\'eductible de $W_{F}$ de dimension $d_{\rho}$, on note $Jord_{\rho}(\psi_{0})$ l'ensemble des entiers $a$ (s'il en existe) tel que $\rho\otimes \sigma_{a}$ soit une sous-repr\'esentation irr\'eductible de $\psi_{0}$.  Si $Jord_{\rho}(\psi_{0})\neq \emptyset$, on note alors  $a_{\rho,\psi_{0}}$ son \'el\'ement maximal.
 Supposons maintenant que $Jord_{\rho}(\psi_{0})=\emptyset$; si  $\rho \otimes \sigma_{2}$ est \`a valeurs dans un groupe de m\^eme type que $G^*$ (cf. ci-dessus); on pose $a_{\rho,\psi_{0}}$ et sinon c'est $\rho$ qui est \`a valeurs dans un groupe de m\^eme type que $G^*$ et on pose $a_{\rho,\psi_{0}}:=-1$. 

On voit maintenant $\rho$ comme une repr\'esentation cuspidale irr\'eductible de $GL(d_{\rho},F)$.
On d\'efinit $x_{\rho,\pi_{0}}$ comme l'unique d'apr\`es \cite{silberger} r\'eel positif ou nul tel que l'induite $\rho
\vert\,\vert^{x_{\rho,\pi_{0}}}\times \pi_{0}$ soit r\'eductible et on montre que (cf. \ref{pointsdereductibilite}, ci-dessous) que $x_{\rho,\pi_{0}}=(a_{\rho,\psi_{0}}+1)/2$.

Au passage, on d\'emontre que $\psi_{0}$ est sans trou c'est-\`a-dire que si $\psi_{0}$ contient une repr\'esentation de la forme $\rho\otimes \sigma_{a}$ avec $a>2$, alors il contient aussi la repr\'esentation $\rho\otimes \sigma_{a-2}$ (ici $\rho$ est une repr\'esentation irr\'eductible de $W_{F}$). Pour 
d\'emontrer le r\'esultat annonc\'e calculant les points de r\'eductibilit\'e, on d\'emontre d'abord que $x_{\rho,\pi_{0}}$ est un demi-entier, puis on v\'erifie tous les cas sauf le dernier cas ci-dessus, en utilisant la compatibilit\'e des paquets \`a l'induction et la restriction (cela a d\'ej\`a \'et\'e fait en \cite{transfert}). Puis on montre que le dernier cas en r\'esulte ''par d\'efaut''.

Avec ces r\'esultats sur les repr\'esentations cuspidales, \cite{europe} et \cite{ams} associent \`a toute 
s\'erie discr\`ete irr\'eductible un morphisme de $W_{F}\times SL(2,{\mathbb C})$ dans $GL(n^*,{\mathbb C})$, not\'e $\psi_{\pi}$, il n'y a plus qu'\`a d\'emontrer la compatibilit\'e avec les r\'esultats d'Arthur, c'est-\`a-dire que $\pi\in \Pi(\psi_{\pi})$; c'est fait en \ref{compatibilite}. 

\

Dans la suite du travail, on se limite comme expliqu\'e ci-dessus au cas o\`u $G=SO(2n+1,F)$ 
(d\'eploy\'e ou non). On suppose que tous les lemmes fondamentaux pour les groupes endoscopiques de $SO(2m+1,F)$ et pour $\tilde{GL}(2m,F)$ sont v\'erifi\'es pour tout $m$. Il y a une autre propri\'et\'e dont on a besoin qui peut s'exprimer tr\`es facilement ainsi: avec la notation $\Pi(\psi)$ d\'ej\`a employ\'ee, il n'y a qu'une combinaison lin\'eaire (\`a homoth\'etie pr\`es) d'\'el\'ements de $\Pi(\psi)$ qui est stable et toute combinaison lin\'eaire stable de s\'eries discr\`etes et dans l'espace vectoriel engendr\'e par ces \'el\'ements. On discute cette propri\'et\'e en \ref{Icuspstable} car il est vraisemblable qu'elle est dans (ou puisse \^etre d\'emontr\'e par les m\'ethodes de) \cite{arthurnouveau}; mais comme elle n'est pas \'ecrite telle quelle, on la prend comme hypoth\`ese. Tel que nous avons \'ecrit les arguments, nous avons besoin de cette hypoth\`ese pour tous les morphismes $\psi'$ \`a valeurs dans $Sp(2m,{\mathbb C})$ pour tout $m\leq n$. 
Avec les lemmes fondamentaux et cette hypoth\`ese, on montre que l'application d\'efinie par Arthur de $\Pi(\psi)$ dans l'ensemble des caract\`eres du groupe $Cent_{Sp(2n,{\mathbb C})}\psi/Cent(Sp(2n,{\mathbb C}))$ est une bijection; un tel r\'esultat est propre au cas des paquets temp\'er\'es. On obtient alors la preuve de notre conjecture sur la classification des repr\'esentations cuspidales de $SO(2m+1,F)$ (pour tout $m\leq n$), \`a savoir cet ensemble de repr\'esentations cuspidales est en bijection avec l'ensemble des couples, $\psi,\epsilon$, form\'e d'un morphisme $\psi$ sans trou de $W_{F}\times SL(2,{\mathbb C})$ dans $Sp(2m,{\mathbb C})$, d\'efinissant une repr\'esentation semi-simple sans multiplicit\'e de $W_{F}\times SL(2,{\mathbb C})$ et d'un caract\`ere $\epsilon$ du centralisateur de $\psi$ dans $Sp(2n,{\mathbb C})$, trivial sur le centre de $Sp(2n,{\mathbb C})$, qui a les propri\'et\'es suivantes:

$\epsilon$ est altern\'e (ou cuspidale au sens de Lusztig) c'est-\`a-dire: soit une
sous-repr\'esentation irr\'eductible $\rho\otimes \sigma_a$ de $W_{F}\times SL(2,{\mathbb C})$ de dimension $ad_{\rho}$ intervenant dans $\psi$; elle est n\'ecessairement \`a valeurs dans un groupe $Sp(ad_{\rho},{\mathbb C})$ dont on note $Z_{\rho\otimes \sigma_a}$ le centre, sous-groupe \`a 2 \'el\'ements qui est naturellement un sous-groupe du centralisateur de $\psi$. Alors $\epsilon$ est altern\'e si pour $\rho,a$ comme ci-dessus, la restriction de $\epsilon$ \`a $Z_{\rho\otimes \sigma_a}\simeq \{\pm 1\}$ est non trivial si $a=2$ et si $a>2$, cette restriction n'est pas le m\^eme caract\`ere que la restriction de $\epsilon$ \`a $Z_{\rho\otimes \sigma_{a-2}}\simeq \{\pm 1\}$;

On remarque que pour cette conjecture pr\'ecise, on n'a pas trait\'e le cas non d\'eploy\'e qui devrait pourtant \^etre analogue en rempla\c{c}ant la condition sur la restriction du caract\`ere au centre de $Sp(2n,{\mathbb C})$ par son oppos\'e. 

\

C'est la derni\`ere conjecture qui manquait pour pouvoir d\'ecrire compl\`etement les paquets de Langlands de s\'erie discr\`ete suivant \cite{europe} et \cite{ams}; on v\'erifie pour finir que nos constructions sont bien compatibles avec celles d'Arthur. Il s'agit ici de v\'erifier que si $\pi\in \Pi(\psi)$, le caract\`ere associ\'e par Arthur se lit sur le module de Jacquet de $\pi$. On montre donc la 
propri\'et\'e suivante, soit comme ci-dessus $\rho\otimes \sigma_a$ une sous-repr\'esentation 
irr\'eductible de $W_{F}\times SL(2,{\mathbb C})$ incluse dans $\psi$. Supposons d'abord qu'il existe $b\in {\mathbb N}$ avec $b<a$ et $\rho\otimes [b]$ une sous-repr\'esentation de $\psi$; on note alors $a_{-}$ le plus grand entier $b$ v\'erifiant ces 2 conditions. Quand un tel $b$ n'existe pas $a_{-}$ n'est pas d\'efini si $a$ est impair et vaut $0$ si $a$ est pair. On note encore $\rho$ la repr\'esentation cuspidale de $GL(d_{\rho},F)$ correspondant \`a $\rho$ par la correspondance de Langlands locale pour 
$GL(d_{\rho},F)$ (\cite{harris},\cite{henniart}). Soit $\pi\in \Pi(\psi)$ et $\epsilon_{\cal A}(\pi)$ le 
caract\`ere du centralisateur de $\psi$ associ\'e par Arthur \`a $\pi$; on montre, avec les notations 
d\'ej\`a introduites ci-dessus, que la restriction de $\epsilon_{\cal A}(\pi)$ \`a $Z_{\rho\otimes \sigma_a}$ est identique   \`a la restriction de ce caract\`ere \`a $Z_{\rho\otimes \sigma_{a_{-}}}$ ( est triviale si $a_{-}=0$) si et seulement si il existe une repr\'esentation $\pi'$ du groupe $SO(2(m-d_{\rho}(a-a_{-})/2)+1,F)$ et une inclusion:
$$
\pi\hookrightarrow \rho\vert\,\vert^{(a-1)/2}\times \cdots \times \rho\vert\,\vert^{(a_{-}+1)/2}\times \pi'.
$$

Je remercie Laurent Clozel pour ces explications sur la globalisation de situations locales et Jean-Loup Waldspurger dont les travaux sont indispensables \`a cet article.

\

L'essentiel des id\'ees encore d\'evelopp\'ees ici et que j'utilise depuis de nombreuses ann\'ees m'ont 
\'et\'e inspir\'ees par l'\'etude de la correspondance de Howe (avec l'interpr\'etation donn\'ee  par J. Adams). C'est donc un immense plaisir pour moi que de pouvoir d\'edier cette article \`a Roger Howe pour son soixanti\`eme anniversaire. Je remercie aussi les organisateurs du congr\`es en son honneur pour leur accueil extr\^emement chaleureux.

\section{Support cuspidal des s\'eries discr\`etes}
\subsection{Notations\label{notations}}
Dans tout l'article sauf la partie \ref{casnondeploye}, si $G(n)$ est un groupe orthogonal, on le suppose d\'eploy\'e.
On reprend les notations de l'introduction, pour tout $m\in {\mathbb N}$, $G(m),G^*(m)$, $m^*$. Et pour tout $m$, on note $J$ la matrice $$\biggl(\begin{matrix} 0 & \cdots &1\\ \vdots& \cdots &\vdots\\
(-1)^{m-1} &\cdots &0
\end{matrix}\bigg)
$$
de fa\c{c}on \`a pouvoir d\'efinir l'action de $\theta$ sur $GL(m,F)$ par $g\mapsto J \, ^tg^{-1}J^{-1}$ pour tout $g\in GL(m,F)$.

\subsection{Quelques rappels des r\'esultats de \cite{arthurnouveau}\label{rappel}}
 On rappelle  des r\'esultats fondamentaux pour nous de \cite{arthurnouveau} 30.1; soit $\psi$ un morphisme continu born\'e de $W_{F}\times SL(2,{\mathbb C})$ dans $G^*(n)$. On le voit comme un morphisme de $W_{F}\times SL(2,{\mathbb C})\times SL(2,{\mathbb C})$ dans $G^*(n)$ trivial sur la deuxi\`eme copie de $SL(2,{\mathbb C})$. Arthur a d\'efini une action de $\theta$ sur la 
 repr\'esentation $\pi(\psi)$ de $GL(n^*,F)$ en imposant \`a $\theta$ d'induire l'action triviale sur l'espace des fonctionnelles de Whittaker (ceci n\'ecessite a priori un choix d'un caract\`ere additif mais on renvoit \`a \cite{transfert} pour l'ind\'ependance et une discussion plus g\'en\'erale); pour donner un sens pr\'ecis \`a cela, il faut fixer, \`a l'aide du caract\`ere addifif fix\'e, un caract\`ere du groupe des matrices unipotentes sup\'erieures invariant par $\theta$; le module des co-invariants de $\pi(\psi)$ pour ce groupe unipotent et ce caract\`ere est de dimension 1 et l'action de $\theta$ est fix\'ee de telle sorte que $\theta$ agisse trivialement dans cet espace de co-invariants. Ceci dit, pour nous ici, la description de la normalisation n'a acune importance. 
 
Arthur en loc. cit. montre alors l'existence d'un paquet fini de repr\'esentations $\Pi(\psi)$ de $G(n)$ et de coefficients entiers positifs $c_{\pi}$ tel que $tr \pi(\psi)\circ \theta$ soit un transfert stable de la distribution $\sum_{\pi\in \Pi(\psi)}c_{\pi}tr\, \pi$; les propri\'et\'es des coefficients sont cach\'es dans la  
d\'efinition de $\tilde{\Pi}_{fin}$ d'Arthur, ils ne joueront pas de r\^ole ici sauf la positivit\'e (\`a un signe commun \`a tous les coefficients pr\`es)! Les 
repr\'esentations incluses dans $\Pi(\psi)$ sont uniquement d\'etermin\'ees gr\^ace \`a l'ind\'ependance 
lin\'eaire des caract\`eres.
On dit que le morphisme $\psi$ est $\theta$-discret si vu comme une repr\'esentation de $W_{F}\times SL(2,{\mathbb C})$ dans $GL(n^*,{\mathbb C})$, il d\'efinit une repr\'esentation sans multiplicit\'e de ce groupe; c'est exactement la m\^eme d\'efinition que celle de l'introduction.

\subsection{Appartenance \`a un paquet stable\label{appartenance}}
Soit $\pi$ une s\'erie discr\`ete irr\'eductible de $G(n)$. Il faut savoir que la projection du caract\`ere distribution de $\pi$ sur l'ensemble des distributions stables \`a support dans les \'el\'ements elliptiques est non nulle. C'est un r\'esultat local et on va en donner 2 d\'emonstrations.

\subsubsection{argument global}
La premi\`ere d\'emonstration repose sur la globalisation suivante qui m'a \'et\'e indiqu\'ee par Laurent Clozel et qui se trouve dans \cite{clozel} theorem 1 B, $\pi$ peut se globaliser, c'est-\`a-dire \^etre consid\'er\'ee comme une composante locale d'une forme automorphe de carr\'e int\'egrable 
$\Sigma$ d'un groupe ad\'elique associ\'e \`a $G$; on peut imposer \`a $\Sigma$ ce que l'on veut en un nombre fini de places, par exemple d'\^etre cuspidale en une autre place et Steinberg en au moins 2 autres places. Cela assure que $\Sigma$ est cuspidale et qu'elle intervient dans une formule des traces simple automatiquement stable; on peut alors d'apr\`es \cite{arthurnouveau} 30.2 (b) trouver une forme automorphe $\tilde{\Sigma}$ du groupe ad\'elique associ\'e \`a $GL(n^*)$ tel que 
$\pi$ soit dans le paquet associ\'e par Arthur \`a la composante locale de $\tilde{\Sigma}$ pour notre place; mais Arthur assure m\^eme qu'en toute place la composante de $\Sigma$ est dans un paquet stable dont un transfert est donn\'ee par la composante en cette m\^eme place de $\Sigma$. On a suppos\'e que $\Sigma $ est Steinberg en au moins une place mais une telle repr\'esentation est automatiquement stable et son seul transfert possible parmi les composantes locales de forme automorphe de carr\'e int\'egrable est la repr\'esentation de Steinberg; on peut le 
v\'erifier avec des arguments de module de Jacquet, ici on insiste sur le fait que le transfert se d\'efinit par  des \'egalit\'es sur toutes les classes de conjugaisons stables pas seulement les elliptiques. Ainsi  $\tilde{\Sigma}$ est n\'ecessairement une repr\'esentation automorphe cuspidal. Et  la composante locale de $\tilde{\Sigma}$ en la place qui nous int\'eresse, est  une repr\'esentation not\'ee $\tilde{\pi}$ de $GL(n^*,F)$  telle que $\pi\in \Pi(\tilde{\pi})$; on sait que $\tilde{\pi}$ est $\theta$-stable mais m\^eme plus qu'elle provient du groupe endoscopique $G(n)$, c'est un des r\'esultats de \cite{arthurnouveau} 30.2.  On ne conna\^{\i}t pas la conjecture de Ramanujan, on \'ecrit donc $\tilde{\pi}$ comme une induite de la  forme 
$$
\times _{(\rho,a,x)\in {\cal I}}St(\rho,a)\vert\,\vert^x\times \tilde{\pi}_{0} \times St(\rho,a)\vert\,\vert^{-x},\eqno(1)
$$ 
o\`u ${\cal I}$ est un ensemble d'indices param\'etrant des triplets form\'es d'une repr\'esen\-tation cuspidale unitaire irr\'educ\-tible $\rho$ d'un $GL(d_{\rho},F)$, d'un entier $a$ et d'un r\'eel $x\in ]0,1/2[$ et o\`u $\tilde{\pi}_{0}$ est temp\'er\'ee et son param\`etre est \`a valeurs (\`a conjugaison pr\`es) dans un groupe $G^*(n_{0})$ pour $n_{0}$ convenable. On pose $\sigma:=\times _{(\rho,a,x)\in {\cal I}}St(\rho,a)\vert\,\vert^x$ et la $\theta$-stabilit\'e entra\^{\i}ne que $$\theta(\sigma)=\times _{(\rho,a,x)\in {\cal I}}St(\rho,a)\vert\,\vert^{-x}.$$
 A $\tilde{\pi_{0}}$, on associe en suivant Arthur 30.1 un paquet de repr\'esentations $\Pi(\tilde{\pi}_{0})$ de $G(n_{0})$ o\`u $n_{0}$ est convenable comme ci-dessus. Pour des coefficients $c_{\pi_{0}}\in {\mathbb R}$ positifs  on a une \'egalit\'e de $$\sum_{\pi_{0}\in \Pi_{0}}c_{\pi_{0}}tr\, \pi_{0}(h')=tr \tilde{\pi}_{0}(g',\theta)\eqno(2)$$
pour tout $h'$ suffisamment r\'egulier dans $G(n_{0})$, o\`u $g',\theta$ a une classe de conjugaison stable qui correspond, dans le transfert, \`a la classe de conjugaison stable de $h'$ (cf \cite{waldspurger1} I.3, III.1, III.2). Les formules explicites pour le calcul de la trace d'une induite donne ici avec les notations de (1):
$$
tr(\sigma \times \tilde{\pi}_{0} \times \theta(\sigma))(g,\theta)= \sum_{\pi_{0}\in \Pi(\tilde{\pi}_{0})}c_{\pi_{0}}tr (\sigma\times \pi_{0}) (h)\eqno(3)
$$
quand la classe de conjugaison stable de $h$ et celle de $(g,\theta)$ se correspondent. Ainsi la repr\'esentation $\pi$ est un sous-quotient irr\'eductible de l'une des induites $\sigma\times \pi_{0}$ pour $\pi_{0}\in \Pi(\pi_{0})$. 
Remarquons aussi pour la suite que si $\psi_{0}$ n'est pas $\theta$ discret alors il existe $\psi'_{0}$ une repr\'esentation de $W_{F}\times SL(2,{\mathbb C})$ dans $GL(m_{0}^*,{\mathbb C})$ se factorisant par $G^*(m_{0})$, $m_{0}$ un entier convenable et une repr\'esentation irr\'eductible $\rho\otimes \sigma_{a}$ de $W_{F}\times SL(2,{\mathbb C})$ tels que $$
\psi_{0}=\psi'_{0}\oplus \bigl(\rho\otimes \sigma_{a} \bigr)\oplus  \bigl(\theta(\rho) \otimes \sigma_{a}\bigr).
$$
Et, comme ci-dessus, on peut \'ecrire de fa\c{c}on imag\'ee 
$$
\Pi(\psi_{0})= St(\rho,a)\times \Pi(\psi'_{0}),
$$
ce qui veut dire que les \'el\'ements du membre de gauche sont tous les sous-quotients irr\'eductibles du membre de droite, c'est-\`a-dire des induites $St(\rho,a)\times \pi'_{0}$ avec $\pi'_{0}\in \Pi(\psi'_{0})$. Par construction les \'el\'ements de $\Pi(\psi'_{0})$ sont des repr\'esentations unitaires et les induites ci-dessus sont donc semi-simples. D\`es que l'on aura d\'emontr\'e que $\sigma$ est n\'ecessairement inexistant, on saura  que, puisque $\pi$ est une s\'erie discr\`ete, $\psi_{0}$ est $\theta$-discret.

 En particulier si $\pi$ est cuspidal, alors $\psi=\psi_{0}$ car $\pi$ ne peut \^etre un sous-quotient de $\sigma\otimes \pi_{0}$ avec $\pi_{0}\in \Pi(\psi_{0})$ et donc $\psi_{0}$ est $\theta$-discret.

 Le d\'efaut de cet argument est qu'il  utilise de fa\c{c}on assez forte  la partie globale des r\'esultats d'Arthur et il n'est donc pas clair \`a priori, qu'il ne faille pas les lemmes fondamentaux pond\'er\'es. On va donc donn\'e un argument local.
 
 \subsubsection{argument d'apparence locale \label{argumentlocal}}
 Ci-dessous tous les arguments sont locaux mais ils utilisent des r\'esultats qui ont \'et\'e obtenu avec des formules des traces simples, donc des arguments globaux. Ils m'ont \'et\'e donn\'es par Waldspurger.

 On sait que le caract\`ere distribution d'une s\'erie discr\`ete, $\pi$, de $G(n)$, peut se voir comme une distribution sur l'ensemble des classes de conjugaison d'\'el\'ements de $G(n)$ dont la restriction aux classes de conjugaison d'\'el\'ements elliptiques est non nulle. En \'evaluant sur ces int\'egrales orbitales, on d\'efinit une application de l'espace vectoriel engendr\'e par ces caract\`eres de repr\'esentations dans un espace $I_{cusp}(G(n))$, l'espace des int\'egrales orbitales des pseudo-coefficients, \'etudi\'e par Arthur en \cite{selecta}; en \'elargissant l'espace des caract\`eres \`a l'ensemble des caract\`eres des repr\'esentations elliptiques de $G(n)$ (d\'efinition d'Arthur) on obtient ainsi une application bijective.
Arthur dans \cite{selecta} a donn\'e une d\'ecomposition en somme directe de cet espace $I_{cusp}(G(n))$, l'un des facteurs \'etant l'espace des distributions stables \`a support dans les \'el\'ements elliptiques.  On a d\'ej\`a v\'efi\'e en \cite{inventiones} que pour d\'efinir la projection sur ce facteur stable, seul le lemme fondamentale stable est n\'ecessaire. On utilise la notation $I^{st}_{cusp}(G(n))$ pour cet espace, conforme aux notations de \cite{waldspurger1} VI. 1.

 Et le point est de d\'emontrer que la projection du caract\`ere distribution de $\pi$ sur cet espace est non nulle; on le fait en localisant au voisinage de l'origine ce qui permet d'utiliser le d\'eveloppement asymptotique d'Harish-Chandra d'o\`u une  d\'ecomposition en somme de transform\'ee de Fourier 
 d'int\'egrales orbitales unipotentes; la dimension des orbites donnent un degr\'e d'homog\'en\'eit\'e pour chacun de ces coefficients. La stabilisation respecte le degr\'e d'homog\'en\'eit\'e. On regarde le coefficient relatif \`a l'orbite triviale, c'est le degr\'e formel de la s\'erie discr\`ete et ce coefficient est donc non nul; la transform\'ee de Fourier de l'orbite de l'\'el\'ement $0$ est automatiquement stable et a donc une projection non nulle sur l'analogue de $I^{st}_{cusp}(G(n))$ dans l'alg\`ebre de Lie. Ceci est le terme de degr\'e 0 dans le developpement asymptotique de la projection du caract\`ere de $\pi$ sur $I^{st}_{cusp}(G(n))$ d'o\`u la non nullit\'e de cette projection.

 On utilise ensuite \cite{waldspurger1} VI.1 proposition (b) qui d\'emontre que le transfert induit un isomorphisme d'espace vectoriel entre $I^{st}_{cusp}(G(n))$ et son analogue stable pour $GL(n^*,F)\times \theta$; cet analogue est 
 d\'efini en \cite{waldspurger1} V.1 o\`u l'existence de ce facteur direct est d\'emontr\'ee. On sait 
 d'apr\`es \cite{waldspurger1} IV.5 (1) que ce dernier espace s'interpr\`ete comme combinaison lin\'eaire de caract\`eres de repr\'esentations temp\'er\'ees $\theta$-stable de $GL(n^*,F)$. Toutefois, \cite{waldspurger1} ne sait pas que les param\`etres de ces
  repr\'esentations sont \`a valeurs dans $G^*(n)$. Pour s'en sortir il faut avoir une description de $I_{cusp}^{st}(G)$; cela est fortement sugg\'er\'e par \cite{arthurnouveau} paragraphe 30 mais n'est pas totalement explicite (cf. notre discussion en \ref{Icuspstable}). Cette m\'ethode est certainement la plus conceptuelle mais se heurte \`a cette difficult\'e. C'est pour cela que l'on a pris l'autre m\'ethode qui 
n\'ecessite de devoir montrer que $\tilde{\pi}$ est temp\'er\'ee ce qui sera fait ci-dessous.
  
   On garde les notations $\psi, \sigma, \psi_{0}$ introduites ci-dessus.

\subsection{Support cuspidal, rappel des d\'efinitions et r\'esultats \'el\'ementaires \label{definitionsupportcuspidal}}

Soit $\pi$ une s\'erie discr\`ete irr\'eductible de $G(n)$. On sait d\'efinir le support cuspidal de $\pi$ (comme de toute repr\'esentation irr\'eductible de $G(n)$) comme l'union d'une repr\'esentation cuspidale irr\'eductible, $\pi_{cusp}$ d'un groupe $G(n_{cusp})$ (le support cuspidal partiel de $\pi$) et d'un ensemble de repr\'esentations cuspidales irr\'eductibles $\rho_{i}$ de groupe $GL(d_{\rho_{i}},F)$ pour $i$ parcourant un ensemble convenable, ${\cal I}$, d'indices. Cet ensemble est d\'efini \`a permutation pr\`es et \`a inversion pr\`es, c'est-\`a-dire que l'on peut changer $\rho_{i}$ en sa duale. Cet ensemble est uniquement d\'efini par la propri\'et\'e que $\pi$ est un sous-quotient irr\'eductible de l'induite 
$$\times_{i\in {\cal I}}\rho_{i}\times \pi_{cusp}.$$

On note $Supp_{GL}(\pi):=\{(\rho'',x)\}$ l'ensemble des couples $\rho'',x$ form\'es d'une 
repr\'esentation cuspidale unitaire d'un groupe lin\'eaire et d'un r\'eel positif ou nul tel que le support cuspidal de $\pi$ soit l'union de $\pi_{cusp}$ avec l'ensemble des repr\'esentations cuspidales $\rho''\vert\,\vert^x$ pour $(\rho'',x)\in Supp_{GL}(\pi)$.

\

\bf Lemme. \sl (i) Soit $\rho''$ une repr\'esentation cuspidale unitaire d'un groupe $GL(d_{\rho''},F)$; on suppose qu'il existe $x''\in {\mathbb R}$ tel que $(\rho'',x'')\in Supp_{GL}(\pi)$. Alors, il existe des r\'eels $d,f$ tels que $d-f+1\in {\mathbb N}_{\geq 1}$ et une s\'erie discr\`ete $\pi'$ tels que $\pi$ soit un sous-module irr\'eductible de l'induite $$<\rho''\vert\,\vert^d, \cdots, \rho''\vert\,\vert^{f}>\times \pi',$$
o\`u la repr\'esentation entre crochet est $St(\rho'',d-f+1)\vert\,\vert^{^{(d+f)/2}}$. De plus $$Supp_{GL}(\pi)=Supp_{GL}(\pi')\cup \{(\rho,\vert y\vert); y\in [d,f]\}.$$

(ii) Pour tout $(\rho'',x)\in Supp_{GL}(\pi)$, $\rho''$ est autoduale et il existe un entier relatif $z$ tel que l'induite $\rho''\vert\,\vert^{x+z}\times \pi_{cusp}$ soit r\'eductible.
\rm

\

Ces r\'esultats ne sont pas nouveaux et bien connus des sp\'ecialistes. Par d\'efinition du support cuspidal pour tout $(\rho',x')\in Supp_{GL}(\pi)$, il existe un choix de signe $\zeta_{(\rho',x')}$ et un ordre sur $Supp_{GL}(\pi)$ tel que l'on ait une inclusion:
$$
\pi \hookrightarrow \times _{(\rho',x')\in Supp_{GL}(\pi)}\rho'\vert\,\vert^{\zeta_{\rho',x'}x'} \times \pi_{cusp}\eqno(1)
$$
On fixe $\rho''$ comme dans l'\'enonc\'e et on prend pour $f$ le plus petit r\'eel de la forme $\zeta_{(\rho'',x'')}x''$ tel que $\rho''\vert\,\vert^f$ interviennent dans (1). On peut ''pousser'' $\rho''\vert\,\vert^f$ vers la gauche tout en gardant une inclusion comme dans (1). Et quand il est le plus \`a gauche possible, toutes les repr\'esentations intervenant avant $\rho''\vert\,\vert^{f}$ sont de la forme $\rho''\vert\,\vert^z$ avec $z$ parcourant un segment de la forme $[d,f[$ pour $d$ convenable v\'erifiant $d-f+1\in {\mathbb N}_{\geq 1}$. Cela entra\^{\i}ne qu'il existe une repr\'esentation $\pi'$ irr\'eductible et une inclusion de $\pi$ dans l'induite $$<\rho''\vert\,\vert^d, \cdots, \rho''\vert\,\vert^{f}>\times \pi'.$$
Il reste \`a d\'emontrer que $\pi'$ est n\'ecessairement une s\'erie discr\`ete. On sait que $d+f>0$ 
car s'il n'en est pas ainsi,  $[d,f]$ est un segment dont le milieu, $(d+f)/2$ est inf\'erieur ou \'egal \`a 
$0$; on pose $\delta=0$ si $(d-f+1)/2$ est un entier et $\delta=(d+f)/2$ sinon. Ainsi
$$
\sum_{y\in [d,f]}y= \sum_{t\in [1,[(d-f+1)/2]} (d-t+1+f+t-1)+\delta= [(d-f+1)/2] (d+f) +(d+f)/2 \leq 0.
$$
Or le module de Jacquet de $\pi$ contient la repr\'esentation $\otimes_{\ell \in [d,f]}\rho\vert\,\vert^\ell \otimes \pi'$ et le crit\`ere de Casselman montre que les exposants sont une combinaison lin\'eaire \`a coefficient strictement positif de toutes les racines simples. L'exposant que l'on vient de trouver ne 
v\'erifie pas cette positivit\'e d'o\`u notre assertion. 
Montrons maintenant que $\pi'$ est une s\'erie discr\`ete; il faut v\'erifier le crit\`ere de Casselman sur les exposants. On suppose a contrario qu'il existe un terme dans le module de Jacquet de $\pi'$ de la forme $\otimes \rho'''\vert\,\vert^{x'''} \otimes \pi_{cusp}$ qui ne satisfait pas au crit\`ere de positivit\'e.  Un tel terme  donne lieu \`a une inclusion de $\pi'$ dans l'induite 
$\times \rho'''\vert\,\vert^{x'''}\times \pi_{cusp}$. Le crit\`ere de positivit\'e de Casselman s'exprime exactement par le fait que pour tout 
$\rho''',x'''$ intervenant dans cette \'ecriture, $\sum_{(\rho,x)} x>0$ o\`u la somme porte sur les 
$(\rho,x)$ \`a gauche de $(\rho''',x''')$ y compris $(\rho''',x''')$. Supposons que $\pi'$ ne soit pas une 
s\'erie discr\`ete. On peut fixer un \'el\'ement de son module de Jacquet ne v\'erifiant pas ce crit\`ere. Et on peut consid\'erer dans $\sum_{(\rho,x)}x $ comme pr\'ec\'edemment, chaque somme partielle o\`u $\rho$ est fix\'e et l'une au moins est $\leq 0$. Par l'argument d\'ej\`a donn\'ee on trouve une inclusion de $\pi'$ dans une induite de la forme:
$$
<\rho'\vert\,\vert^{d'}, \cdots, \rho'\vert\,\vert^{f'}>\times \pi'',
$$
o\`u $[d',f']$ est un segment (i.e. $d'-f'+1\in {\mathbb N}_{\geq 1}$), $\rho'$ est une repr\'esentation cuspidale unitaire et $\sum_{x\in [d',f']}x \leq 0$. Comme ci-dessus, ceci est \'equivalent \`a $d'+f'\leq 0$. L'induite $$<\rho''\vert\,\vert^d, \cdots, \rho''\vert\,\vert^{f}>\times <\rho'\vert\,\vert^{d'}, \cdots, \rho'\vert\,\vert^{f'}>
$$est n\'ecessairement irr\'edutible: en effet soit $\rho''\not\simeq \rho'$ et le r\'esultat est imm\'ediat. Soit $\rho''=\rho'$, mais alors $f\leq f' $ par minimalit\'e de $f$, et 
$$
d-d'=(d+f)-(d'+f')+(f'-f)>0,
$$
d'o\`u $[d,f]\supset [d',f']$ et le r\'esultat. On obtient donc une inclusion de $\pi$ dans l'induite:
$$<\rho''\vert\,\vert^d, \cdots, \rho''\vert\,\vert^{f}>\times
<\rho'\vert\,\vert^{d'}, \cdots, \rho'\vert\,\vert^{f'}>\times \pi'' \simeq
$$
$$
<\rho'\vert\,\vert^{d'}, \cdots, \rho'\vert\,\vert^{f'}>\times <\rho''\vert\,\vert^d, \cdots, \rho''\vert\,\vert^{f}>\times
\pi''
$$
et $d'+f'\leq 0$  contredit le fait que $\pi$ est une s\'erie discr\`ete. L'assertion sur le $Supp_{GL}(\pi)$ est claire par d\'efinition.

(ii) se d\'eduit de (i) par r\'ecurrence: on reprend les notations de (i) et on admet que tous les \'el\'ements de $Supp_{GL}(\pi')$ sont de la forme $(\rho,x)$ avec $\rho$ autoduale. Si $\rho''$ n'est pas autoduale, $\rho''\vert\,\vert^x\times \pi_{cusp}$ est irr\'eductible pour tout $x$ r\'eel par un r\'esultat 
g\'en\'eral d'Harish-Chandra (cf. \cite{JIMJ}) et $\rho''\vert\,\vert^x$ commute avec toute repr\'esentation de la forme $\rho'\vert\,\vert^{\pm x'}$ pour $(\rho',x')\in Supp_{GL}(\pi')$. On aurait donc un isomorphisme:
$$
<\rho''\vert\,\vert^d,\cdots, \rho''\vert\,\vert^f>\times \pi' \simeq <(\rho'')^{*-f}, \cdots, (\rho'')^{*-d}>\times \pi'.
$$
et une inclusion de $\pi$ dans l'induite de droite. Or $-f-d<0$ ce qui contredit le fait que $\pi$ est une s\'erie discr\`ete et le crit\`ere de Casselman. C'est exactement le m\^eme argument pour la 
deuxi\`eme partie de (ii)

\subsection{Propri\'et\'es de demi-integralit\'e du support cuspidal d'une s\'erie discr\`ete\label{integralite}}
\bf Th\'eor\`eme. \sl (i) Soit $\pi$ une s\'erie discr\`ete irr\'eductible de $G(n)$; tout \'el\'ement $(\rho,x)$ de $Supp_{GL}(\pi)$ (notation de \ref{definitionsupportcuspidal}) est tel que $\rho$ est autoduale et que $x$ est demi-entier.

(ii) Soit $\pi_{0}$ une repr\'esentation cuspidale de $G(n)$ et $\rho$ une repr\'esentation cuspidale 
irr\'eductible autoduale de $GL(d_{\rho},F)$ et $x\in {\mathbb R}$ tel que l'induite $\rho\vert\,\vert^x\times \pi_{0}$ soit r\'eductible. Alors $x\in 1/2 {\mathbb Z}$.

\rm

\

Gr\^ace \`a \ref{definitionsupportcuspidal} (ii), pour d\'emontrer (i),  il suffit de d\'emontrer que pour $\pi_{cusp}$ une repr\'esentation cuspidale irr\'eductible d'un groupe $G(n_{cusp})$ et pour $\rho$ une 
repr\'esentation cuspidale irr\'eductible d'un groupe $GL(d_{\rho},F)$, autoduale, le r\'eel positif ou nul, $x_{\rho,\pi_{cusp}}$ tel que l'induite $\rho\vert\,\vert^{x_{\rho,\pi_{cusp}}}\times \pi_{cusp}$ est 
r\'eductible, est  un demi-entier. C'est-\`a-dire, en fait, (ii) et c'est donc (ii) que nous allons d\'emontrer.

Soit donc $\pi_{cusp}$ et $x_{\rho,\pi_{cusp}}$ comme ci-dessus. Si $x_{\rho,\pi_{cusp}}=0$, l'assertion est claire. Supposons donc que $x_{\rho,\pi_{cusp}}>0$. Alors l'induite:
$$
\rho\vert\,\vert^{x_{\rho,\pi_{cusp}}+1}\times \rho\vert\,\vert^{x_{\rho,\pi_{cusp}}}\times \pi_{cusp}
$$
 a un  sous-module irr\'eductible qui est une s\'erie discr\`ete; c'est un calcul de module de Jacquet facile on coince cette sous-repr\'esentation dans l'intersection des 2 sous-repr\'esentations suivantes: (la notation $<\sigma_{1},\sigma_{2}>$ repr\'esente le socle de l'induite $\sigma_{1}\times \sigma_{2}$, c'est \`a dire la somme des sous-modules irr\'eductibles, mais ici le socle est irr\'eductible)
 $$
 \rho\vert\,\vert^{x_{\rho,\pi_{cusp}}+1}\times <\rho\vert\,\vert^{x_{\rho,\pi_{cusp}}},\pi_{cusp}> \cap
 <\rho\vert\,\vert^{x_{\rho,\pi_{cusp}}+1},\rho\vert\,\vert^{x_{\rho,\pi_{cusp}}}>\times \pi_{cusp}.
 $$
 Puis on calcule les modules de Jacquet, calcul qui prouve d'abord que cette intersection est non nulle puis que c'est une s\'erie discr\`ete.
 
 On note $\pi$ cette s\'erie discr\`ete et on lui applique \ref{appartenance} avec les notations de ce paragraphe. En particulier, si $\psi\neq \psi_{0}$,  $Supp_{GL}(\pi)$ contient le support cuspidal de 
 $\sigma$ (\`a ''conjugaison'' pr\`es). Mais $Supp_{GL}(\pi)$ a exactement 2 \'el\'ements $\rho\vert\,\vert^{x_{\rho,\pi_{cusp}}+1}$ et $\rho\vert\,\vert^{x_{\rho,\pi_{cusp}}}$; cela entra\^{\i}ne
 d\'ej\`a que $\sum_{(\rho,a,x)\in {\cal I}}a\leq 2$. Si ${\cal I}$ est r\'eduit \`a un \'el\'ement $(\rho,a,x)$ avec $a=2$,  les ensembles $((1+x_{\rho,\pi_{cusp}}),x_{\rho,\pi_{cusp}})$ et $((1/2+x), (-1/2+x))$ 
 co\"{\i}ncident \`a l'ordre et au signe pr\`es. D'o\`u par positivit\'e, puisque $x<1/2$
 $$
 1+x_{\rho,\pi_{cusp}}=1/2+x; \quad x_{\rho,\pi_{cusp}}=1/2-x.
 $$
 Ceci est impossible. On ne pas  non plus avoir  $\vert {\cal I}\vert =2$. On sait donc maintenant que soit $\psi=\psi_{0}$ soit  $\sigma$ est de la forme $\rho\vert\,\vert^{x}$; on rappelle que $x\in ]0,1/2[$ et qu'il co\"{\i}ncide donc avec $x_{\rho,\pi_{cusp}}$ et non avec $1+x_{\rho,\pi_{cusp}}$. En d'autres termes soit $\pi\in \Pi(\tilde{\pi})$ avec $\tilde{\pi}=\tilde{\pi}_{0}$ est temp\'er\'ee soit un sous-quotient irr\'eductible de l'induite $\rho\vert\,\vert^{1+x_{\rho,\pi_{cusp}}}\times \pi_{cusp}$ est dans $\Pi(\tilde{\pi}_{0})$. En fait l'induite est irr\'eductible par d\'efinition de $x_{\rho,\pi_{cusp}}$ et dans ce dernier cas, ce serait toute l'induite qui serait dans $\Pi(\tilde{\pi}_{0})$.

 On montre maintenant que  tout $\pi_{0}\in \Pi(\tilde{\pi}_{0})$ est tel que $Supp_{GL}(\pi_{0})$ est 
 form\'e de couples $(\rho',x')$ avec $x'$ demi-entier. Cela permettra donc de conclure gr\^ace \`a ce qui pr\'ec\`ede; on aura en plus montr\'e que $\sigma$ est trivial.

 On fixe $\rho$ une repr\'esentation cuspidale unitaire irr\'eductible d'un groupe $GL(d_{\rho},F)$.
 Pour $z\in {\mathbb R}$, $m\in {\mathbb N}$ et $\pi'$ une repr\'esentation de $G(m)$ on note $Jac_{z}\pi'$ l'unique \'el\'ement du groupe de Grothendieck des repr\'esentations lisses de longueur finie de $G(m)$ tel que la restriction de $\pi'$ au Levi $GL(d_{\rho},F)\times G(m)$ le long du radical unipotent d'un parabolique de ce Levi soit de la forme $\rho\vert\,\vert^z\otimes Jac_{z}\pi' \oplus \tau$ o\`u $\tau$ est une somme de repr\'esentations irr\'eductibles de la forme $\rho'\otimes \tau''$ avec $\rho'\not\simeq \rho\vert\,\vert^z$. Pour $\tilde{\pi}$ une repr\'esentation de $GL(m^*+2d_{\rho'},F)$ on d\'efinit $Jac^\theta_{z}\tilde{\pi}$ comme l'unique \'el\'ement du groupe de Grothendieck des 
repr\'esentations lisses de longueur finie de $GL(m,F)$ tel que la restriction de $\tilde{\pi}$ au Levi $GL(d_{\rho'},F) \times GL(m^*,F) \times GL(d_{\rho'},F)$ soit de la forme 
$\rho\vert\,\vert^z\otimes Jac^\theta_{z}\tau \otimes \rho^*\vert\,\vert^{-z}\oplus \tau'$ o\`u $\tau'$ est une somme de repr\'esentations irr\'eductibles de la forme $\rho'\otimes \tau'' \otimes \rho''$ avec soit $\rho'\not\simeq \rho\vert\,\vert^z$ soit $\rho''\not\simeq \rho^{-z}$. On a vu en \cite{inventiones} (mais c'est assez facile) que l'\'egalit\'e des traces \ref{appartenance} (1) donne aussi une \'egalit\'e (avec les m\^emes coefficients):
$$
\sum_{\pi\in \Pi(\tilde{\pi}_{0})}c_{\pi_{0}} tr Jac_{z}\pi_{0} (h)=tr Jac^\theta_{z} \tilde{\pi}_{0}(g,\theta),\eqno(1)
$$
Le terme de droite ne peut \^etre non nul que pour des $z\in 1/2 {\mathbb Z}$.
Ici on utilise la positivit\'e des $c_{\pi_{0}}$ annonc\'ee par Arthur qui emp\^eche toute simplification dans le terme de gauche. Alors  pour $\pi_{0}\in \Pi(\tilde{\pi}_{0})$, $Jac_{z}\pi_{0}\neq 0$ n\'ecessite que $Jac^\theta_{z}\tilde{\pi}_{0}\neq 0$. D'o\`u le fait que $z\in 1/2{\mathbb Z}$. Si $Supp_{GL}(\pi_{0})$ contient  un \'el\'ement $(\rho'',x)$ avec $x\in {\mathbb R}$ ,de fa\c{c}on standard (en poussant vers la gauche), on montre qu'il existe $x'$ avec $x-x'\in {\mathbb Z}$ tel que $Jac_{x'}\pi_{0}\neq 0$, en faisant ici $\rho=\rho'$. Comme on vient de montrer que $x'\in 1/2 {\mathbb Z}$, $x$ aussi est un demi-entier relatif.

Cela termine la preuve du th\'eor\`eme.  
 \subsection{Morphisme associ\'e \`a une s\'erie discr\`ete\label{morphismeassocie}}
Le corollaire suivant est \'evidemment tr\`es fortement inspir\'e par les travaux d'Arthur.

\bf Corollaire. \sl Soit $\pi$ une s\'erie discr\`ete irr\'eductible de $G(n)$, alors il existe une 
repr\'esentation temp\'er\'ee, $\tilde{\pi}$ de $GL(n^*,F)$, $\theta$-discr\`ete  telle que $\pi\in \Pi(\tilde{\pi})$. Ou encore, il existe un morphisme $\psi$ de $W_{F}\times SL(2,{\mathbb C})$ dans $GL(n^*,{\mathbb C})$ tel que $\pi\in \Pi(\psi)$ avec $\psi$ un morphisme $\theta$-discret \`a valeurs dans $G^*(n)$.
\rm

En tenant compte de \ref{appartenance} il faut montrer que $\tilde{\pi}=\tilde{\pi}_{0}$ avec les notations de ce paragraphe. Comme ci-dessus, n\'ecessairement $Supp_{GL}(\pi)$ contient le support cuspidal de $\sigma$ (\`a ''conjugaison'' pr\`es). Mais maintenant on sait que $Supp_{GL}(\pi)$ est form\'e de couple $(\rho,x)$ avec $x$ demi-entier et ce n'est pas le cas du support cuspidal de $\sigma$ si $\sigma$ existe vraiment. D'o\`u le fait que $\tilde{\pi}=\tilde{\pi}_{0}$. Le corollaire r\'esulte alors de \ref{appartenance}.

\section{Morphismes associ\'es aux repr\'esentations cuspidales de $G(n)$ et points de r\'eductibilit\'e des induites de cuspidales\label{pointsdereductibilite}.}
, 

On dit qu'un morphisme  $\psi$ de $W_{F}\times SL(2,{\mathbb C})$ dans un $GL(m,{\mathbb C})$ est sans trou si pour toute repr\'esentation $\rho\otimes \sigma_{a}$ intervenant dans $\psi$ (notations de l'introduction) avec $a>2$ la repr\'esentation $\rho\otimes \sigma_{a-2}$ y intervient aussi.

On pr\'ecise les notations de l'introduction; soit $\pi_{0}$ une repr\'esentation cuspidale de $G(m)$ et soit $\rho$ une repr\'esentation cuspidale irr\'eductible autoduale de $GL(d_{\rho},F)$ ce qui d\'efinit $d_{\rho}$. On note $x_{\rho,\pi_{0}}$ l'unique r\'eel positif ou nul (cf. \cite{silberger}) tel que l'induite $\rho\vert\,\vert^{x_{\rho,\pi_{0}}}\times \pi_{0}$ soit r\'eductible. Soit $\psi_{0}$ tel que $\pi_{0}\in \Pi(\psi_{0})$ ce qui est possible d'apr\`es le th\'eor\`eme pr\'ec\'edent. Pour $\rho$ comme ci-dessus, on pose 

$Jord_{\rho}(\psi_{0}):=\{a\in {\mathbb N}; $ tel que la repr\'esentation $\rho\otimes \sigma_{a}$  soit une sous-repr\'esentation de $\psi_{0}\}$.

\noindent On note $a_{\rho,\psi_{0}}$ l'\'el\'ement maximal de $Jord_{\rho}(\psi_{0})$ quand cet ensemble est non vide. On dit qu'une repr\'esentation de $W_{F}\times SL(2,{\mathbb C})$ (ou simplement de $W_{F}$) est du type de $G^*$ si elle est orthogonale quand $V$ est symplectique et symplectique quand $V$ est orthogonal.

\

\bf Corollaire. \sl Soient $\rho,\pi_{0}, \psi_{0}$ comme ci-dessus et tels que $\pi_{0}\in \Pi(\psi_{0})$ (cf. \ref{morphismeassocie}). Alors $\psi_{0}$ est $\theta$ discret et sans trou.

(i)
On suppose que $Jord_{\rho}(\psi_{0})\neq \emptyset$ alors $x_{\rho,\pi_{0}}=(a_{\rho,\psi_{0}}+1)/2$. 

(ii) On suppose que $Jord_{\rho}(\psi_{0})=\emptyset$ mais que $\rho\otimes \sigma_{2}$ est une 
repr\'esentation du type de $G^*$; alors $x_{\rho,\pi_{0}}=1/2$.

(iii) On suppose que $Jord_{\rho}(\psi_{0})=\emptyset$ et que $\rho$ vue comme repr\'esentation de $W_{F}$ est de m\^eme type que $G^*$. Alors $x_{\rho,\pi_{0}}=0$.

(iv) $Jord(\pi_{0})$ (cf introduction pour la notation) co\"{\i}ncide avec $Jord(\psi_{0})$. En particulier $\psi_{0}$ est uniquement d\'etermin\'e par $\pi_{0}$.\rm

\

Le fait que $\psi_{0}$ est $\theta$ discret est dans \ref{appartenance}. 
Le fait que $\psi_{0}$ soit sans trou et les assertions (i) et (ii) sont d\'emontr\'ees en \cite{transfert}. On va red\'emontrer l'essentielle de ces assertions pour la commodit\'e du lecteur;  on ne redonne pas la d\'emonstration de ce que si $Jord_{\rho}(\psi_{0})\neq \emptyset$ alors $x_{\rho,\pi_{0}}=(a_{\rho,\psi_{0}}+1)/2$; la d\'emonstration est plus simple mais de m\^eme nature que celle de (ii) que l'on va redonner.

En \cite {algebra}, cf. l'introduction, on a montr\'e l'in\'egalit\'e, o\`u la somme porte sur toutes les 
repr\'esentations cuspidales irr\'eductibles autoduales d'un $GL(d,F)$:
$$
\sum_{\rho}\sum_{\ell \in [1,[x_{\rho,\pi_{0}]]}}d_{\rho}(2x_{\rho}-2\ell+1)\leq m^*.
$$
Cette somme ne voit que les $\rho$ tel que $x_{\rho}>1/2$. Or si on se limite dans le terme de gauche aux repr\'esentations $\rho$ telles que $Jord_{\rho}(\psi_{0})\neq \emptyset$, ce terme est sup\'erieur ou \'egal \`a $\sum_{(\rho,a)\in Jord(\psi)}ad_{\rho}$ avec \'egalit\'e si et seulement si $Jord(\psi_{0})$ est sans trou. Or cette somme vaut $m^*$car c'est exactement la somme des dimensions des sous-repr\'esentations de $\psi_{0}$. Et on a donc d\'emontr\'e que $\psi_{0}$ est sans trou et que pour tout $\rho$ tel que $Jord_{\rho}(\psi_{0})=\emptyset$, $x_{\rho,\pi_{0}}\in [0,1/2]$. Comme on sait d\'ej\`a que $x_{\rho,\pi_{0}}$ est un demi-entier, il reste les possibilit\'es $0$ ou $1/2$. 
On remarque ici que le fait que $x_{\rho,\pi_{0}}$ soit entier ou non r\'esulte de la parit\'e de $a_{\rho,\psi_{0}}$ et donc de savoir si le morphisme de $W_{F}$ dans $GL(d_{\rho},{\mathbb C})$  associ\'e \`a $\rho$ est du type de $G^*$ ou non. Ce sont les conditions de parit\'es voulues dans la d\'efinition des blocs de Jordan.

Avec cela (iv) r\'esulte de \cite{algebra} 4.3.

 Il reste donc \`a montrer, pour toute repr\'esentation autoduale cuspidale irr\'eductible, $\rho$, telle que $Jord_{\rho}(\psi_{0})=\emptyset$,  l'\'equivalence entre le fait que $\rho$ est telle que $x_{\rho,\pi_{0}}=1/2$ et le fait que $\rho\otimes \sigma_{2}$ est de m\^eme type que $G^*$. On rappelle que $\sigma_{2}$ est une repr\'esentation symplectique et donc que $\rho\otimes \sigma_{2}$ est une repr\'esentation symplectique si $\rho$ est une repr\'esentation orthogonale et est une repr\'esentation orthogonale si $\rho$ est une repr\'esentation symplectique.

Soit donc $\rho$ tel que $x_{\rho,\pi_{0}}=1/2$; on note ici $\pi_{d}$ l'unique sous-module irr\'eductible de l'induite $\rho\vert\,\vert^{1/2}\times \pi_{0}$; c'est une s\'erie discr\`ete par hypoth\`ese de 
r\'eductibilit\'e. On consid\`ere une repr\'esentation temp\'er\'ee, $\theta$-discr\`ete, $\tilde{\pi}_{d}$ telle que $\pi\in \Pi(\tilde{\pi}_{d})$ et on \'ecrit avec des coefficients $c_{\pi'}$ positifs, comme en \ref{appartenance}
$$
\sum_{\pi'\in \Pi(\tilde{\pi}_{d})}c_{\pi'}tr \pi' (h)= tr \tilde{\pi}_{d}(g,\theta).
$$
Comme $\tilde{\pi}_{d}$ est une repr\'esentation temp\'er\'ee, on l'\'ecrit comme une induite  de repr\'esentations de Steinberg g\'en\'eralis\'ees, $\times_{(\rho',a)\in {\cal I}}St(\rho',a)$ o\`u ${\cal I}$ est un ensemble d'indices param\'etrant des couples $(\rho',a)$ form\'es d'une repr\'esentation cuspidale unitaire $\rho'$ et d'un entier $a$; on sait que ${\cal I}$ est sans multiplicit\'e puisque $\tilde{\pi}_{d}$ est $\theta$-discret. Avec les notations $Jac_{z}$ et $Jac^\theta_{z}$ de la preuve de \ref{morphismeassocie} pour $\rho$ fix\'e comme ici, on a encore que $Jac^\theta_{1/2}\tilde{\pi}_{d}$ est un transfert de $$\sum_{\pi'\in \Pi(\tilde{\pi})}c_{\pi'}tr Jac_{1/2}\pi'.$$
Or $Jac_{1/2}\pi_{d}=\pi_{0}$, ainsi $Jac^\theta_{1/2}\tilde{\pi}_{d}\neq 0$. Cela n\'ecessite que ${\cal I}$ contienne $(\rho,2)$; il le contient avec multiplicit\'e exactement 1. On obtient alors:
$$
Jac^\theta_{1/2}\tilde{\pi}_{d}=\times_{(\rho',a)\in {\cal I}-\{(\rho,2)\}}St(\rho',a)=:\tilde{\pi}_{0}
$$
et $\pi_{0}\in \Pi(\tilde{\pi}_{0})$. On note encore $\psi_{0}$ le morphisme associ\'e \`a $\pi_{0}$ et le param\`etre de $\tilde{\pi}_{d}$ est donc la somme de $\psi_{0}$ et de $\rho\otimes \sigma_{2}$ (ici $\rho$ est la repr\'esentation irr\'eductible de $W_{F}$ associ\'ee \`a $\rho$ via la correspondance de Langlands et $\sigma_{2}$ est la repr\'esentation irr\'eductible de dimension 2 de $SL(2,{\mathbb C})$). D'apr\`es la description d'Arthur (formule (30.15) de 30.2) ce paquet doit \^etre de m\^eme type que $G^*$; comme c'est d\'ej\`a le cas de $\psi_{0}$, il en est de m\^eme de $\rho\otimes \sigma_{2}$.

On suppose maintenant que $\rho$ est tel que $\rho\otimes \sigma_{2}$ est de m\^eme type que $G^*$ et que $Jord_{\rho}(\psi)=\emptyset$; on doit montrer que $x_{\rho,\pi_{0}}=1/2$.  On note $\psi_{d}:=\rho\otimes \sigma_{2}\oplus \psi_{0}$ le morphisme de $W_{F}\times SL(2,{\mathbb C})$ dans $GL(m^*+2d_{\rho},{\mathbb C})$ et $\tilde{\pi}_{d}$ la repr\'esentation temp\'er\'ee de $GL(m^*+2d_{\rho},F)$ lui correspondant. Le morphisme $\psi_{d}$ est \`a valeurs dans $G^*(m+d_{\rho})$. On a donc encore un transfert, pour des bons coefficients positifs (\cite{arthurnouveau} 30.1)
$$
\sum_{\pi'\in \Pi(\psi_{d})}c_{\pi'} tr\, \pi' (h)= tr \tilde{\pi}_{d}(g,\theta).\eqno(1)
$$
Et avec les notations $Jac_{z}$ et $Jac^\theta_{z}$ d\'ej\`a introduite pour $\rho$ fix\'e et $z\in {\mathbb R}$:
$$
\sum_{\pi'\in \Pi(\psi_{d})}c_{\pi'}Jac_{z} tr\, \pi'(h')=Jac^\theta_{z} tr \tilde{\pi}_{d}(g',\theta).\eqno(2)
$$
On l'applique d'abord avec $z=1/2$; le terme de droite n'est pas nul et il contient $\tilde{\pi}_{0}$; il est 
m\^eme r\'eduit \`a $\tilde{\pi}_{0}$ car par hypoth\`ese $Jord_{\rho}(\pi_{0})=0$ et on a d\'ej\`a identifi\'e $Jord_{\rho}(\pi_{0})$ avec les entiers $a$ tel que $\rho\otimes \sigma_{a}$ soit une sous-repr\'esentation de $\psi_{0}$. Donc $\psi_{0}$ ne contient pas $\rho\otimes \sigma_{2}$. Ainsi le terme de gauche de (2) est le paquet associ\'e \`a $\psi_{0}$ et il contient donc $\pi_{0}$. Ainsi il existe $\pi'\in \Pi(\psi_{d})$ tel que $Jac_{1/2}\pi'$ contienne $\pi_{0}$ ou encore que $\pi'$ est un sous-quotient irr\'eductible de $\rho\vert\,\vert^{1/2} \times \pi_{0}$. Comme $Jac_{1/2}\pi'\neq 0$, $\pi'$ est un sous-module irr\'eductible de l'induite $\rho\vert\,\vert^{1/2}\times \pi_{0}$. Cette induite est n\'ecessairement 
r\'eductible car sinon $Jac_{-1/2}\pi'\neq 0$ alors que pour $z=-1/2$ le terme de droite de (2) est nul. Ainsi $x_{\rho,\pi_{0}}=1/2$ comme cherch\'e et cela termine la d\'emonstration.

\section{Classification et paquet de Langlands\label{compatibilite}}
Soit $\pi$ une s\'erie discr\`ete irr\'eductible  de $G(n)$; nous lui avons associ\'e un morphisme de 
$W_{F}\times SL(2,{\mathbb C})$ dans $GL(n^*,{\mathbb C})$, \`a l'aide de ses blocs de Jordan. On note $\psi_{\pi}$ ce morphisme. En suivant Arthur, on lui a aussi associ\'e un morphisme en \ref{morphismeassocie}, le th\'eor\`eme ci-dessous dit que ces morphismes sont conjugu\'es:

\

\bf Th\'eor\`eme. \sl On a $\pi\in \Pi(\psi_{\pi})$.\rm

\

Ce th\'eor\`eme montre que nos constructions sont compatibles avec celles d'Arthur et que le morphisme  associ\'e \`a une s\'erie discr\`ete \`a l'aide de ses blocs de Jordan est bien celui conjectur\'e par Langlands. Toutefois, il n'est pas d\'emontr\'e que les coefficients permettant de construire la distribution stable associ\'ee au paquet sont \'egaux \`a 1; ce probl\`eme est r\'esolu par Waldspurger en \cite{waldspurger2} pour certaines repr\'esentations venant des constructions de Lusztig (\cite{lusztig}) et pour $G(n)=SO(2n+1,F)$. Et il n'est pas non plus d\'emontr\'e que les repr\'esentations \`a l'int\'erieur du paquet sont classifi\'ees par les caract\`eres du centralisateur; ce point semble moins 
s\'erieux et nous r\`eglerons le cas des groupes orthogonaux d\'eploy\'es ici avec une m\'ethode tout \`a fait g\'en\'erale.
\

Pour \'eviter les confusions, on note  $\psi$ plut\^ot que $\psi_{\pi}$ le morphisme associ\'e \`a $\pi$ par  \cite{europe} et \cite{ams} \`a l'aide des blocs de Jordan; on rappellera les propri\'et\'es qui le 
caract\'erise ci-dessous.  Et on note $\tilde{\psi}$ celui associ\'e essentiellement par Arthur en \ref{morphismeassocie} et $\tilde{\pi}$ la repr\'esentation temp\'er\'ee de $GL(n^*,F)$ associ\'e \`a $\tilde{\psi}$. Il faut d\'emontrer que $\psi$ et $\tilde{\psi}$ sont conjugu\'es. On r\'eutilise les notations $Jac_{z}$ et $Jac^\theta_{z}$ de \ref{integralite}. On v\'erifie d'abord que $\tilde{\psi}$ est $\theta$-discret. Sinon $\pi$ serait un sous-quotient d'une induite de la forme $St(\rho,a)\times \pi'$ avec $\pi'$ une repr\'esentation au moins unitaire; l'induite est alors semi-simple et $\pi$ est temp\'er\'e mais non discret. Ainsi pour tout $z\in {\mathbb R}$, $Jac^\theta_{z}\tilde{\pi}$ est 0 ou est une repr\'esentation temp\'er\'ee  irr\'eductible.
On a vu en loc. cit. que pour tout $z\in {\mathbb R}$, 
$$
Jac_{z}\pi \in \Pi(Jac^\theta_{z}\tilde{\pi}).
$$

Si $\pi$ est une repr\'esentation cuspidale, le fait que $\psi$ et $\tilde{\psi}$ sont conjugu\'es a \'et\'e vu en \ref{pointsdereductibilite} (iv). On suppose donc que $\pi$ n'est pas cuspidale.
D'apr\`es la construction de \cite{europe},  2 cas sont \`a distinguer. Dans ce qui suit $\rho$ est une repr\'esentation cuspidale autoduale 
irr\'eductible d'un $GL(d_{\rho},F)$.

1e cas: il existe $x_{0}$ de la forme $(a-1)/2$ avec $a\in {\mathbb N}_{>1}$ et une s\'erie discr\`ete $\pi'$ tel que $\pi$ soit l'unique sous-module irr\'eductible de l'induite $\rho\vert\,\vert^{x_{0}}\times \pi'$. Dans ce cas, si on note $\psi'$ le morphisme de $W_{F}\times SL(2,{\mathbb C})$ dans $G^{n^*-d_{\rho}}$, alors $\psi'$ contient $(\rho,a-2)$ comme bloc de Jordan ($a=2$ est accept\'e) et $\psi$ s'obtient en rempla\c{c}ant ce bloc de Jordan par $(\rho,a)$.

Concluons dans ce cas; comme $Jac^\theta_{x=(a-1)/2}\tilde{\pi}_{v}\neq 0$ et est 
irr\'eductible, cela veut dire que la repr\'esentation $\rho\otimes \sigma_{a}$ intervient avec multiplicit\'e exactement 1 dans $\tilde{\psi}$ et que la repr\'esentation obtenue s'obtient simplement en rempla\c{c}ant $\tilde{\psi}$ par un morphisme $\tilde{\psi'}$ ayant les m\^emes sous-repr\'esentations irr\'eductibles sauf cette repr\'esentation qui devient $\rho\otimes \sigma_{a-2}$. Il suffit d'appliquer par exemple une hypoth\`ese de 
r\'ecurrence pour savoir que $\tilde{\psi'}$ et $\psi'$ sont conjugu\'es pour obtenir la 
m\^eme assertion pour $\tilde{\psi}$ et $\psi$.

2e cas: il existe $x_{0}$ de la forme $(a+1)/2$ avec $a\in {\mathbb N}_{\geq 1}$ et une 
s\'erie discr\`ete $\pi'$ tel que $\pi$ soit l'un des 2 sous-modules irr\'eductibles de l'induite:
$$
<\rho\vert\,\vert^{(a+1)/2},St(\rho,a)>\times \pi'.
$$
On note $\psi'$ le morphisme associ\'e \`a $\pi'$ par nos construction et on admet encore par 
r\'ecurrence que $\pi'\in \Pi(\psi')$. La diff\'erence avec le cas 1 est que $Jac_{x_{0}}\pi$ est une repr\'esentation $\pi_{1}$ qui est une repr\'esentation temp\'er\'ee et non plus une s\'erie discr\`ete. Toutefois, nous avons d\'ej\`a v\'erifi\'e que puisque $\pi'\in \Pi(\psi')$, alors $\pi_{1}\in \Pi(\psi'\oplus (\rho\otimes \sigma_{a} \oplus \rho\otimes \sigma_{a}))$ et comme ci-dessus, $$\tilde{\psi}=\psi'\oplus \rho\otimes \sigma_{a+2}\oplus \rho\otimes \sigma_{a}.$$ Mais d'apr\`es nos constructions, on a bien que $Jord(\psi)$ se d\'eduit de $Jord(\psi')$ en ajoutant les 2 blocs $(\rho,a)$ et $(\rho,a+2)$ ce qui est exactement la 
m\^eme chose. Ceci prouve le th\'eor\`eme.

\section{Classification \`a la Langlands des s\'eries discr\`etes de $SO(2n+1,F)$.}
Dans cette partie on suppose que  $G(n)$ est la forme d\'eploy\'ee du groupe $SO(2n+1,F)$. On a en vue la classification de Langlands des s\'eries discr\`etes. Soit $\psi$ un morphisme $\theta$-discret de $W_{F}\times SL(2,{\mathbb C})$ dans $GL(2n,{\mathbb C})$, on dit qu'il est $\theta$-stable s'il se factorise par $Sp(2n,{\mathbb C})$. 
\subsection{D\'efinition de l'espace $I_{cusp}(\,)$ \label{definitiondeIcusp}}
On note $C_{cusp}(G)$ l'ensemble des fonctions lisses cuspidales de $G$, c'est-\`a-dire celles dont les int\'egrales orbitales sur les \'el\'ements non elliptiques sont nulles et $I_{cusp}(G)$ est l'image de cet espace vectoriel de fonctions modulo le sous-espace des fonctions lisses dont toutes les int\'egrales orbitales sont nulles.
L'article d'Arthur \cite{selecta} montre toute l'importance de cet espace; dans cet article le groupe est suppos\'e connexe, ce qui est le cas de $G$.  On d\'efinit $I_{cusp}^{st}(G)$ comme l'ensemble des 
\'el\'ements qui ont des int\'egrales orbitales constantes sur les classes de conjugaison stable. On appelle $I_{cusp}^{ns}(G)$ l'ensemble des \'el\'ements dont la somme des int\'egrales orbitales sur toute classe de conjugaison stable est nulle. Waldspurger a remarqu\'e (cf. \cite{inventiones} 4.5) que sans aucun lemme fondamental, on peut montrer la d\'ecomposition:
$$
I_{cusp}(G)=I_{cusp}^{st}(G)\oplus I_{cusp}^{ns}(G).\eqno(1)
$$
Le point est que la stabilit\'e commute \`a la transformation de Fourier \cite{fourier}. Pour aller plus loin,
on suppose la validit\'e des lemmes fondamentaux pour les groupes endoscopiques de $G$  les lemmes fondamentaux et pour les alg\`ebres de Lie, introduits par Waldspurger. Pour toute donn\'ee endoscopique, not\'e abusivement $H$ de $G$, on peut d'une part d\'efinir $I_{cusp}^{st}(H)$ (en rempla\c{c}ant $G$ par $H$ dans les notations ci-dessus) et un transfert de $I_{cusp}^{st}(H)$ dans $I_{cusp}(G)$; on note $I_{cusp}^{H-st}(G)$ l'image; cette image peut se d\'efinir directement. Et le r\'esultat principal de \cite{selecta} est de prouver l'\'egalit\'e:
$$
I_{cusp}^{ns}(G)=\oplus_{H}I_{cusp}^{H-st}(G).\eqno(2)
$$
Waldspurger d\'efinit et \'etudie l'espace analogue \`a $I_{cusp}(G)$ pour les groupes non connexes de la forme $\tilde{GL}(m',F)$ et montre en \cite{waldspurger1} VI.1 que sous l'hypoth\`ese d'un lemme fondamental convenable, le transfert induit un isomorphisme de $I_{cusp}^{st}(G)$ sur $I_{cusp}^{st}(\tilde{GL}(2n,F))$. On peut ainsi r\'ecrire (1) et (2) en
$$
I_{cusp}(G)=\oplus_{H}I_{cusp}^{H-st}(G), \eqno(3)
$$
o\`u ici $H$ parcourt toutes les donn\'ees endoscopiques en incluant pour la partie stable $H=\tilde{GL}(2n,F)$.

\subsection{Repr\'esentations elliptiques et $I_{cusp}$ \label{Icusppsi}}
Soit $\pi$ une s\'erie discr\`ete; on sait qu'elle poss\`ede un pseudo-coefficient et on peut d\'efinir la projection de ce pseudo coefficient sur l'ensemble des fonctions cuspidales; pour le cas le plus nouveau celui de $\tilde{GL}(2m,F)$ ceci est expliqu\'e en \cite{waldspurger1} II. Ceci s'\'etend aux 
repr\'esentations elliptiques et permet de d\'efinir pour toute repr\'esentation elliptique un \'el\'ement de $I_{cusp}(G)$; $I_{cusp}(G)$ est engendr\'e comme espace vectoriel par ces images.

Soit maintenant $\psi$ un morphisme de $W_{F}\times SL(2,{\mathbb C})$ dans $Sp(2n,{\mathbb C})$; on supppose que $\psi$ est $\theta$-discret, c'est-\`a-dire que la repr\'esentation d\'efinie quand on inclut $Sp(2n,{\mathbb C})$ dans $GL(2n,{\mathbb C})$ est sans multiplicit\'e; pour aller plus vite on appelle un tel morphisme $\theta$-stable (car il est \`a image dans $Sp(2n,{\mathbb C})$) et $\theta$-discret. On a alors v\'erifi\'e que les \'el\'ements de $\Pi(\psi)$ sont des s\'eries discr\`etes. 

On suppose maintenant que  $\psi$ n'est pas   $\theta$-discret; $\Pi(\psi)$ est  toujours d\'efini uniquement par la propri\'et\'e de transfert stable. Puisque $\psi$ n'est pas $\theta$-discret, $\psi$ est \`a valeurs dans un sous-groupe de Levi d'un parabolique $\theta$-stable et les \'el\'ements de $\Pi(\psi)$ sont les sous-modules irr\'eductibles d'une induite de repr\'esentations temp\'er\'ees; en particulier $\Pi(\psi)$ ne contient pas de s\'eries discr\`etes.

En g\'en\'eral, on note $I_{cusp}(G)[\psi ]$ l'espace vectoriel engendr\'e par l'image dans $I_{cusp}(G)$ des 
repr\'esentations elliptiques combinaisons lin\'eaires de repr\'esentations incluses dans $\Pi(\psi)$; cet espace peut-\^etre 0. Mais si $\psi$ est $\theta$-discret  on a: $$
dim\, I_{cusp}(G)[\psi]= \vert \Pi(\psi)\vert. \eqno(1)
$$
On a v\'erifi\'e en \ref{compatibilite} que pour $\pi$ une s\'erie discr\`ete, un morphisme $\psi$ tel que $\pi\in \Pi(\psi)$ est uniquement d\'etermin\'e par les blocs de Jordan de $\pi$; en particulier si $\Pi(\psi)$ contient une s\'erie discr\`ete alors pour tout morphisme $\psi'$ $\theta$-stable, $\Pi(\psi)\cap \Pi(\psi')=\emptyset$.

On pose $I_{cusp}(G)[nd]$ la somme des espaces $I_{cusp}(G)[\psi]$ pour tous les morphismes $\psi$ non $\theta$-discret. On a donc la d\'ecomposition en somme directe:
$$
I_{cusp}(G)=\oplus_{\psi}I_{cusp}(G)[\psi ] \oplus I_{cusp}(G)[nd]. \eqno(2)
$$
\subsection{Description de $I_{cusp}^{st}(G)$\label{Icuspstable}}
Pour tout morphisme $\psi$ $\theta$-stable et $\theta$-discret, en \cite{arthurnouveau}
30.1 (30.14) (pour $s=1$), il est construit un \'el\'ement de $I_{cusp}^{st}(G)[\psi]$. On a besoin d'avoir
que $I_{cusp}^{st}(G)$ est engendr\'e par ces \'el\'ements. Ceci n'est pas \'ecrit dans \cite{arthurnouveau} mais est tr\`es fortement sugg\'er\'e par les r\'esultats de loc.cit. pour les raisons suivantes.

Pour tout morphisme $\psi$, $\theta$-discret mais non n\'ecessairement $\theta$-stable, on note $I_{cusp}(\tilde{G})[\psi]$ le sous-espace vectoriel de dimension 1 de $I_{cusp}(\tilde{GL}(2n,F))$ engendr\'e par l'image d'un prolongement \`a $\tilde{GL}(2n,F)$ de la repr\'esentation temp\'er\'ee de $GL(2n,F)$ d\'efinie par $\psi$. Waldspurger a montr\'e que $I_{cusp}(\tilde{GL}(2n,F))$ est 
engendr\'e par ces sous-espaces vectoriels. Le probl\`eme est donc de savoir si $I_{cusp}^{st}(\tilde{GL}(2n,F))$ est engendr\'e par ceux de ces \'el\'ements qui correspondent aux $\psi$, $\theta$-stable. A $\psi$ simplement $\theta$-discret, Arthur lui-m\^eme associe un groupe endoscopique elliptique de $\tilde{GL}(2m,F)$ en \cite{arthurnouveau} pages 235 et 236: on d\'ecompose $\psi$ en la somme de 2 morphismes $\psi_{s}\oplus \psi_{o}$ o\`u $\psi_{s}$ est \`a valeurs dans un groupe symplectique, $Sp(2m_{s},{\mathbb C})$ alors que $\psi_{o}$ est \`a valeurs dans un groupe orthogonal, $O(2m_{o},{\mathbb C})$. 
Le d\'eterminant de la restriction de $\psi_{o}$ \`a $W_{F}$ donne un caract\`ere quadratique $\eta$ du corps de base (par r\'eciprocit\'e) et le groupe endoscopique est d\'etermin\'e par $m_{s}, m_{o}$ et $\eta$. On note $H_{\psi}$ ce groupe endoscopique, c'est le groupe $Sp(2m_{s},F)\times O(2m_{o},F)$ o\`u le groupe orthogonal est le groupe de la  forme orthogonale de dimension $2m_{o}$ de discriminant (normalis\'e) $\eta$ et d'invariant de Hasse est $+1$. Par 30.1, Arthur associe \`a $\psi$ un \'el\'ement stable de $I_{cusp}(H_{\psi})$ puisque l'on est dans une situation produit mais il ne transfert pas cet \'el\'ement en un 
\'el\'ement de $I_{cusp}^{H_{\psi}-st}(\tilde{GL}(2m,F))$. Comme on admet les lemmes fondamentaux on peut bien faire ce transfert et le point est de montrer que ce transfert co\"{\i}ncide avec l'image de $\pi(\psi)\circ \theta$ dans $I_{cusp}(\tilde{G}_{m})$.
On admet donc pour la suite de cette partie l'hypoth\`ese suivante:

\

\sl $I_{cusp}^{st}(\tilde{GL}(2n,F))$ est la somme des espaces $I_{cusp}(\tilde{G})[\psi ]$, o\`u $\psi$ parcourt l'ensemble des morphismes $\theta$-stables et $\theta$-discrets.\rm

\

La cons\'equence de cette hypoth\`ese est que la d\'ecomposition \ref{Icusppsi} (2) est compatible \`a la projection sur $I_{cusp}^{st}(G)$ et pr\'ecis\'ement que l'on a:
$$
I_{cusp}^{st}(G)=\oplus_{\psi}\biggl( I_{cusp}(G)[\psi ] \cap I_{cusp}^{st}(G)\biggr), \eqno(1)
$$
o\`u $\psi$ parcourt l'ensemble des classes de conjugaison de morphismes $\theta$-stables et $\theta$-discrets.

\

Montrons que cette hypoth\`ese entra\^{\i}ne aussi la d\'ecomposition en somme directe pour tout $\psi$ morphisme $\theta$-stable et $\theta$-discret:

$$
I_{cusp}(G)[\psi]=\oplus_{H} \biggl(I_{cusp}^{H-st}(G)\cap I_{cusp}(G)[\psi]\biggr), \eqno(2)
$$o\`u $H$ parcourt le m\^eme ensemble qu'en \ref{definitiondeIcusp} (3).

Ceci n'\'etant pas directement dans \cite{arthurnouveau} 30.1, il faut le v\'erifier. Et il suffit de prouver l'analogue de (1) pour $H$ un groupe endoscopique de $Sp(2m,F)$; pour un tel groupe il faut aussi tenir compte de $I_{cusp}(G)[nd]$.
Soit $H$  un produit de groupes orthogonaux de la forme $SO(2m_{1}+1,F)\times SO(2m_{2}+1,F)$ avec $m=m_{1}+m_{2}$; on applique la propri\'et\'e (1) \`a chaque facteur pour d\'ecrire $I_{cusp}^{st}(H)$ comme somme directe des espaces vectoriels de dimension 1, $I_{cusp}^{st}(H)[\psi_{1}\times \psi_{2}]$, o\`u $\psi_{i}$, pour $i=1,2$ est un morphisme $\theta$-stable et $\theta$-discret \`a valeurs dans $Sp(2m_{i},{\mathbb C})$. On sait transf\'erer un tel espace vectoriel dans $I_{cusp}(G)$ gr\^ace \`a \cite{arthurnouveau} 30.1 (30.14); le transfert est \`a valeurs dans $I_{cusp}(G)(\psi_{1,2})$ o\`u $\psi_{1,2}$ est le morphisme de $W_{F}\times SL(2,{\mathbb C})$ dans $Sp(2m,{\mathbb C})$ obtenu en composant $\psi_{1}\times \psi_{2}$ avec l'inclusion naturelle de $Sp(2m_{1},{\mathbb C})\times Sp(2m_{2},{\mathbb C})$ dans $Sp(2m,{\mathbb C})$. Cela donne exactement la d\'ecomposition:
$$
I_{cusp}^{H-st}(G)=\oplus_{\psi}\biggl(I_{cusp}(G)[\psi ]\cap I_{cusp}^{H-st}(G)\biggr) \oplus \biggl(I_{cusp}(G)[nd]\cap I_{cusp}^{H-st}(G)\biggr).
$$
On en d\'eduit la double somme:
$$
I_{cusp}(G)=\oplus_{H} \biggl(\oplus_{\psi}(I_{cusp}(G)[\psi]\cap I_{cusp}^{H-st}(G))\oplus I_{cusp}(G)[nd]\cap I_{cusp}^{H-st}(G)\biggr).
$$
Et on peut inverser les sommations pour obtenir le r\'esultat cherch\'e (2).

\subsection{Classification des s\'eries discr\`etes\label{classification}}
On rappelle que l'on a admis tous les lemmes fondamentaux de notre situation et l'hypoth\`ese de \ref{Icuspstable}.

\

\bf Th\'eor\`eme. \sl Soit $\psi$ un morphisme $\theta$-stable et $\theta$-discret que l'on voit comme une repr\'esentation de $W_{F}\times SL(2,{\mathbb C})$ \`a valeurs dans $GL(2m,{\mathbb C})$. On note $\ell_{\psi}$ le nombre de sous-repr\'esentations irr\'eductibles incluses dans $\psi$. Alors:
$$
\vert \Pi(\psi)\vert=2^{\ell_{\psi}-1}.
$$Et l'application d\'efinie par Arthur qui associe un caract\`ere du centralisateur de $\psi$ dans $Sp(2n,{\mathbb C})$ est une bijection sur l'ensemble des caract\`eres de restriction triviale au centre de $Sp(2n,{\mathbb C})$.

\

\rm
On sait d\'ej\`a que $I_{cusp}^{st}(G)[\psi ]$ est de dimension 1. On fixe $H=SO(2m_{1}+1,F)\times SO(2m_{2}+1,F)$ et on doit calculer le nombre de couples $\psi_{1},\psi_{2}$ tels que $\psi_{i}$ est $\theta$-stable et $\theta$-discret pour le groupe $SO(2m_{i}+1,F)$ (pour $i=1,2$) et tels $\psi$ soit conjugu\'e de $\psi_{1}\times \psi_{2}$ puisque chacun de ces couples donnent un \'el\'ement de 
$I_{cusp}^{H-st}(G)[\psi]$ et que cela l'engendre compl\`etement.  Il est plus simple de faire ce calcul en laissant varier la d\'ecomposition de $m$ en $m_{1}+m_{2}$; ensuite il faut diviser le r\'esultat par 2 \`a cause des isomorphismes entre groupes endoscopiques. Ici on autorise $m_{1}m_{2}=0$, 
c'est-\`a-dire que l'on retrouve 2 fois la partie correspondant \`a $I_{cusp}^{st}(G)[\psi]$. On note $Jord(\psi)$ l'ensemble des sous-repr\'esentations irr\'eductibles de $W_{F}\times SL(2,{\mathbb C})$ incluses dans $\psi$; on rappelle que c'est un ensemble sans multiplicit\'e. Le nombre de d\'ecomposition de $\psi$ en $\psi_{1}\times \psi_{2}$ est pr\'ecis\'ement le nombre de d\'ecomposition de $Jord(\psi)$ en 2 sous-ensembles, c'est \`a dire $2^{\ell_{\psi}}$ puisque $\ell_{\psi}$ est le nombre d'\'el\'ements de $Jord(\psi)$ par 
d\'efinition. On obtient le nombre d'\'el\'ement de $\Pi(\psi)$ quand on a divis\'e par 2.

On v\'erifie ais\'ement que le cardinal du centralisateur de $\psi$ dans $Sp(2m,{\mathbb C})$ est 
$2^{\ell_{\psi}}$ et que le centre est de cardinal 2. Ainsi le groupe des caract\`eres du centralisateur de $\psi$ dans $Sp(2m,{\mathbb C})$ triviaux sur le centre de $Sp(2m,{\mathbb C})$ est aussi de cardinal $2^{\ell_{\psi}-1}$.
Pour d\'emontrer la deuxi\`eme assertion du th\'eor\`eme, il suffit donc de montrer que l'application 
d\'efinie par Arthur est surjective. On note $\epsilon_{\cal A}(\pi)$ cette application. Ce que l'on 
conna\^{\i}t est le rang de la matrice dont les lignes sont ind\'ex\'ees par les \'el\'ements de $Centr_{Sp(2m,{\mathbb C})}\psi/Cent (Sp(2m,{\mathbb C}))$ et les colonnes par les \'el\'ements de $\Pi(\psi)$, les coefficients de la matrice \'etant $$\epsilon_{\cal A}(\pi)(s); 
s\in Centr_{Sp(2m,{\mathbb C})}\psi/Cent (Sp(2m,{\mathbb C})),\pi\in \Pi(\psi).$$
C'est une matrice carr\'e.
V\'erifions que le rang de cette matrice est le nombre de ses lignes: fixons une ligne donc $s \in Centr_{Sp(2m,{\mathbb C})}\psi/Cent (Sp(2m,{\mathbb C}))$. Un tel \'el\'ement fixe une 
d\'ecomposition de $\psi$ en $\psi_{1}\times \psi_{2}$, o\`u pour $i=1,2$, $\psi_{i}$ est inclus 
dans l'espace propre pour l'une des 2 valeurs propres de $s$; comme $s$ n'est d\'efini que modulo le centre de $Sp(2m,{\mathbb C})$, cette d\'ecomposition n'est d\'efinie qu'\`a l'\'echange pr\`es des 2 facteurs. L'\'el\'ement $s$ d\'efinit aussi un groupe endoscopique, $H_{s}$, \`a isomorphisme 
pr\`es (qui est le groupe $\tilde{GL}(2m,F)$ si $s$ est central). Et l'\'el\'ement $$\sum_{\pi} \epsilon_{\cal A}(\pi)(s) \, \pi\eqno(*)_{s}$$ a pour image dans $I_{cusp}(G)$ l'\'el\'ement \cite{arthurnouveau} (30.14) correspondant \`a $H_{s}$ et \`a la d\'ecomposition $\psi=\psi_{1}\times \psi_{2}$ de $\psi$. 
R\'eciproquement \`a toute d\'ecomposition de $\psi$ en $\psi_{1}\times \psi_{2}$, on associe un \'el\'ement $s$ simplement en donnant ses 2 espaces propres comme ci-dessus. Ainsi les \'el\'ements 
$(*)_{s}$ donne la base de  $I_{cusp}(G)[\psi]$ d\'ej\`a consid\'er\'ee. Ils sont lin\'eairement 
ind\'ependants et forment un ensemble de cardinal $2^{\ell_{\psi}-1}$ puisque ce nombre est la dimension de l'espace vectoriel $I_{cusp}(G)[\psi ]$. Cela prouve notre assertion sur le rang de la matrice et termine la preuve.

\subsection{Caract\`ere et module de Jacquet\label{caractereetmoduledejacquet}}
On fixe $\psi$ un morphisme $\theta$-stable et $\theta$-discret.
On reprend les notations de l'introduction; soit $\rho\otimes \sigma_a$ une sous-repr\'esentation irr\'eductible de $\psi$ vu comme repr\'esentation de $W_{F}\times SL(2,{\mathbb C})$. Dans l'introduction on a 
d\'efini l'entier $a_{-}$ dans les 2 cas suivants:

soit il existe $b\in {\mathbb N}$ tel que $\rho\otimes \sigma_{b}$ soit une sous-repr\'esentation irr\'eductible de $\psi$ avec $b<a_{-}$, auquel cas on pose $a_{-}$ le plus grand des entiers $b$ avec ces propri\'et\'es

soit il n'existe pas d'entier $b$ comme ci-dessus mais $a$ est pair et on pose $a_{-}=0$.

On a not\'e $Z_{\rho\otimes \sigma_a}$ le sous-groupe \`a 2 \'el\'ements du centralisateur de $\psi$ dans $Sp(2n,{\mathbb C})$ qui correspond \`a cette sous-repr\'esentation $\rho\otimes \sigma_a$; pr\'ecis\'ement cette repr\'esentation est \`a valeurs dans un sous-groupe symplectique de $Sp(2n,{\mathbb C})$ et $Z_{\rho\otimes \sigma_a}$ est le centre de ce sous-groupe.

Soit $\pi\in \Pi(\psi)$ et notons $\epsilon_{\cal A}(\pi)$ le caract\`ere du centralisateur de $\psi$ qu'Arthur associe \`a $\pi$; $\epsilon_{\cal A}(\pi)$ est connu quand on conna\^{\i}t toutes ses restrictions aux sous-groupes $Z_{\rho\otimes \sigma_a}$. On avait remarqu\'e en \cite{europe} que ce caract\`ere devait \^etre li\'e aux propri\'et\'es de modules de Jacquet de $\pi$. C'\'etait le point de d\'epart des classifications de \cite{europe} et \cite{ams} et on va montrer que $\epsilon_{\cal A}(\pi)$ a bien la propri\'et\'e qui permet de 
d\'efinir les caract\`eres en loc. cit. ou encore que $\epsilon_{\cal A}(\pi)$ co\"{\i}ncide avec le caract\`ere associ\'e \`a $\pi$ par \cite{europe} l\`a o\`u nous l'avions d\'efini.

Pour unifier le th\'eor\`eme ci-dessous, on dit, par convention que la restriction d'un caract\`ere \`a $Z_{\rho\otimes \sigma_{a_{-}}}$ est le caract\`ere trivial si $a_{-}=0$, cas o\`u $Z_{\rho\otimes \sigma_{a_{-}}}$ n'est pas le groupe $\{\pm 1\}$. On note $d_{\rho}$ la dimension de la repr\'esentation $\rho$ et on note encore $\rho$ la repr\'esentation cuspidale de $GL(d_{\rho},F)$ associ\'ee \`a $\rho$ par la correspondance de Langlands.

\

\bf Th\'eor\`eme. \sl Soit $\rho\otimes \sigma_a$ une sous-repr\'esentation de $\psi$ telle que $a_{-}$ soit 
d\'efini. Alors, la restriction de $\epsilon_{\cal A}(\pi)$ \`a $Z_{\rho\otimes \sigma_a}\simeq \{ \pm 1\}$ est le 
m\^eme caract\`ere que la restriction de $\epsilon_{\cal A}(\pi)$ \`a $Z_{\rho\otimes \sigma_{a_{-}}}$ si et seulement si il existe une repr\'esentation $\pi'$ du groupe $SO(2n-d_{\rho}(a+a_{-})+1,F)$ et une inclusion 
$$
\pi\hookrightarrow \rho\vert\,\vert^{(a-1)/2}\times \cdots \times \rho\vert\,\vert^{(a_{-}+1)/2}\times \pi'.
$$
\rm
La cons\'equence de ce th\'eor\`eme est que l'on conna\^{\i}t tr\`es explicitement les modules de Jacquet des repr\'esentations dans $\Pi(\psi)$. 

\

On  consid\`ere la donn\'ee endoscopique de $G$ dont le groupe $H$ est $SO(d_{\rho}(a+a_{-})+1,F)$ $\times$ $SO(2n-d_{\rho}(a+a_{-})+1,F)$. Et pour $H$ on consid\`ere le morphisme $\psi_{1}\times \psi_{2}$ de $W_{F}\times SL(2,{\mathbb C})$ \`a valeurs dans $Sp(d_{\rho}(a+a_{-}),{\mathbb C}) \times Sp(2n-d_{\rho}(a+a_{-}),{\mathbb C})$, o\`u $\psi_{1}$ est la somme $\rho\otimes \sigma_a\oplus \rho\otimes \sigma_{a_{-}}$ et $\psi_{2}$  est la somme des autres sous-repr\'esentations incluses dans $\psi$.

On a besoin de conna\^{\i}tre le paquet $\Pi(\psi_{1})$; on sait par \ref{classification} qu'il a 2 \'el\'ements si $a_{-}\neq 0$ et 1 \'el\'ement sinon. La situation est donc particuli\`erement simple. On sait que 
l'\'el\'ement de $I_{cusp}^{st}(\tilde{GL}(d_{\rho}(a+a_{-}),F))$ qui d\'efinit ce paquet est l'image d'un prolongement \`a $\tilde{GL}(d_{\rho}(a+a_{-}),F)$ de la repr\'esentation temp\'er\'ee $St(\rho,a)\times St(\rho,a_{-})$. On \'ecrit $St(\rho,a)\times St(\rho,a_{-})$ comme l'unique sous-module irr\'eductible de l'induite:
$$
\rho\vert\,\vert^{(a-1)/2}\times \cdots \times \rho\vert\,\vert^{(a_{-}+1)/2}\times St(\rho,a_{-})\times 
St(\rho,a_{-}) \times \rho\vert\,\vert^{-(a_{-}+1)/2} \times \cdots \times \rho\vert\,\vert^{-(a-1)/2}.
$$
Alors $\Pi(\psi_{1})$ contient les sous-modules irr\'eductibles de l'induite pour $SO(d_{\rho}(a+a_{-})+1,F)$:
$$
\rho\vert\,\vert^{(a-1)/2}\times \cdots \times \rho\vert\,\vert^{(a_{-}+1)/2}\times St(\rho,a_{-}).\eqno(1)
$$
C'est un calcul de module de Jacquet expliqu\'e en \cite{transfert} 4.2 et analogue \`a ceux fait ici; on montre en loc.cit. qu'un module de Jacquet convenable de la distribution stable dans $\Pi(\psi_{1})$ (relativement au facteur $GL(d_{\rho},F)$) est non nul et a pour transfert stable $St(\rho,a_{-})\times St(\rho,a_{-})$; avec les notations de loc.cit. c'est $Jac_{(a-1)/2, \cdots, (a_{-}+1)/2}$. Ce module de Jacquet est donc l'induite \`a $SO(2a_{-}+1,F)$ de la repr\'esentation $St(\rho,a_{-})$.  Si $a_{-}>0$, cette induite est de longueur 2 car la parit\'e de $a_{-}$ est la bonne (cf. \ref{pointsdereductibilite}) et $\Pi(\psi_{1})$ contient n\'ecessairement les 2 sous-modules de (1). 
Si $a_{-}=0$, (1) a un unique sous-module irr\'eductible n\'ecessairement dans $\Pi(\psi_{1})$.  On a donc d\'ecrit $\Pi(\psi_{1})$ comme l'ensemble des sous-modules irr\'eductibles de (1) et quand il y en a 2, la distribution stable est la somme de ces 2 sous-modules puisqu'un module de Jacquet d'une distribution stable est stable (cf. \cite{transfert}4.2)

On ne peut \'evidemment pas donner une description aussi pr\'ecise de $\Pi(\psi_{2})$ et de la distribution stable qui est form\'ee avec ses \'el\'ements; on fixe donc des nombres complexes $c_{\pi_{2}}$ pour tout $\pi_{2}\in \Pi(\psi_{2})$ tel que $\sum_{\pi_{2}\in \Pi(\psi_{2})}c_{\pi_2}\pi_{2}$ soit stable. Ainsi la distribution stable correspondant \`a $\psi_{1}\times \psi_{2}$ est:
$$
(\sum_{\pi_{1}\in \Pi(\psi_{1})}\pi_{1})\otimes (\sum_{\pi_{2}\in \Pi(\psi_{2})}c_{\pi_{2}}tr\, \pi_{2}).\eqno(2)
$$
Pour tout $\pi\in \Pi(\psi)$, on fixe $c_{\pi}\in {\mathbb C}$ tel que $\sum_{\pi\in \Pi(\psi)}c_{\pi}\pi$ soit stable. Et on pose $\epsilon_{a}(\pi):=\epsilon_{\cal A}(\pi)(z_{\rho\otimes \sigma_a})$, o\`u $z_{\rho\otimes \sigma_a}$ est l'\'el\'ement non trivial de $Z_{\rho\otimes \sigma_a}$; on d\'efinit de m\^eme $\epsilon_{a_{-}}(\pi)$ si $a_{-}\neq 0$ sinon on pose $\epsilon_{a_{-}}(\pi)=1$.
D'apr\`es \cite{arthurnouveau} (30.14), la  distribution   (2) se 
transf\`ere en une distribution $H$-stable de $G$ de la forme (\`a un scalaire pr\`es qui vient de 
l'impr\'ecision dans le choix des $c_{\pi}$)
$$
\sum_{\pi\in \Pi(\psi)}\epsilon_{a}(\pi)\epsilon_{a_{-}}(\pi)c_{\pi}tr\, \pi.\eqno(3)
$$
Ecrivons explicitement ce transfert sur les caract\`eres vus comme fonctions localement $L^1$ sur les \'el\'ements semi-simples. Pour tout $\gamma_{G}$ \'el\'ement semi-simple de $G$, on l'\'egalit\'e o\`u les $\Delta(\gamma_{G},\gamma_{H})$ sont les facteurs de transfert et o\`u $\gamma_{H}$ parcourt un ensemble de repr\'esentant de classes de conjugaison dans $H$ dont la classe stable se 
transf\`ere en celle de $\gamma_{G}$:
$$
\biggl(\sum_{\pi\in \Pi(\psi)}\epsilon_{a}(\pi)\epsilon_{a_{-}}(\pi) c_{\pi}tr\, \pi\biggr)(\gamma_{G})=
$$
$$
\sum_{\gamma_{H}}\Delta(\gamma_{G},\gamma_{H})\biggl((\sum_{\pi_{1}\in \Pi(\psi_{1})}tr\, \pi_{1})(\sum_{\pi_{2}\in \Pi(\psi_{2})}c_{\pi_{2}}tr\, \pi_{2})\biggr)(\gamma_{H}).\eqno(4)
$$
On va appliquer cette \'egalit\'e en imposant \`a $\gamma_{G}$ d'\^etre dans le sous-groupe de Levi $M$ $\simeq$ $GL(d_{\rho},F) \times SO(2(n-d_{\rho})+1,F)$ de $G$. On \'ecrit un tel \'el\'ement sous la forme $\gamma_{G}=m\times \gamma'$, o\`u $m\in GL(d_{\rho},F)$. 

On note $M^1_{H}$ le sous-groupe de Levi de $H$ isomorphe \`a
$$
GL(d_{\rho},F)\times SO(2(n_{1}-d_{\rho})+1,F) \times SO(2n_{2},F)$$et $M^2_{H}$ celui qui est isomorphe \`a 
$$
SO(2n_{1}+1,F)\times GL(d_{\rho},F)\times SO(2(n_{2}-d_{\rho})+1,F).
$$
Fixons $\gamma_{G}\in M$ comme ci-dessus mais on suppose que $m$ est elliptique dans 
$GL(d_{\rho},F)$. Dans (4) il suffit alors de sommer sur les \'el\'ements $\gamma_{H}$ avec les propri\'et\'es impos\'ees qui sont soit dans $M^1_{H}$ soit dans $M^2_{H}$. On note $Z_{M}$ le centre de $M$ et on applique cela non pas \`a $\gamma_{G}$ mais \`a l'ensemble des $\gamma_{G}z$ o\`u $z$ parcourt $Z_{M}$. On prenant une limite convenable et on appliquant une formule de Casselman 
on peut remplacer dans (4), $tr\, \pi$ par $tr\, res_{M}(\pi)$ et les produits $$\sum_{\gamma_{H}}\Delta(\gamma_{G},\gamma_{H})(tr\, \pi_{1} tr\, \pi_{2})(\gamma_{H})$$ par $$\sum_{\gamma^1_{H}}\Delta (\gamma_{G},\gamma^1_{H})tr\,( res_{M^1_{H}}\pi_{1})tr\, \pi_{2}(\gamma^1_{H})+ \sum_{\gamma^2_{H}} \Delta (\gamma_{G},\gamma^2_{H}) tr\, \pi_{1} tr\, (res_{M_{H}^2}\pi_{2})(\gamma^2_{H}),$$o\`u $\gamma^1_{H}$ est dans $M^1_{H}$ et $\gamma^2_{H}$ est dans $M^2_{H}$.

En faisant varier $m$ dans l'ensemble des \'el\'ements elliptiques de $GL(d_{\rho},F)$ on projette sur le caract\`ere de la repr\'esentation cuspidale $\rho\vert\,\vert^{(a-1)/2}$; parce que $\rho$ est cuspidale bien que l'on se limite aux \'el\'ements elliptiques, cela permet de faire dispara\^{\i}tre toutes les autres repr\'esentations c'est-\`a-dire remplacer les restrictions par les restrictions suivies par cette projection sur $\rho\vert\,\vert^{(a-1)/2}$. Il faut v\'erifier que cette op\'eration appliqu\'ee \`a n'importe quel \'el\'ement de $\Pi(\psi_{2})$ donne $0$. Cela r\'esulte des propri\'et\'es standards 
des modules de Jacquet des paquets de s\'eries discr\`etes que nous avons \'etablies; pour que cette projection soit non nulle il faudrait que la sous-repr\'esentation $\rho\otimes \sigma_a$ de $W_{F}\times SL(2,{\mathbb C})$ soit incluse dans $\psi_{2}$, ce qui n'est pas le cas par choix de $\psi_{2}$. Comme nous allons en avoir besoin ci-dessous remarquons que cette propri\'et\'e 
reste vraie en rempla\c{c}ant $(a-1)/2$ par $(a'-1)/2$ pour tout $a'\in ]a_{-},a]$ car $\rho\otimes \sigma_{a'}$ pour un tel $a'$ n'est pas sous repr\'esentation de $\psi$ par minimalit\'e de $a_{-}$.
Si $a>a_{-}+2$, on recommence ces op\'erations pour le Levi $M'=GL(d_{\rho},F)\times SO(2(n-2d_{\rho})+1,F)$ en projetant cette fois sur $\rho\vert\,\vert^{(a-3)/2}$. On note $$M_{a,a_{-}}\simeq 
GL(d_{\rho},F)\times \cdots \times GL(d_{\rho},F) \times SO(2n-d_{\rho}(a-a_{-})+1,F)$$ le sous-groupe de Levi de $G$ isomorphe o\`u il y a $(a-a_{-})/2$ copies de $GL(d_{\rho},F)$. Et on note
 $$ M^1_{a,a_{-}}= Gl(d_{\rho},F)\times \cdots \times GL(d_{\rho},F)\times SO(2n_{1}-d_{\rho}(a-a_{-})+1,F)\times SO(2n_{2}+1,F)$$
le sous-groupe de Levi de $H$ o\`u il y a $(a-a_{-})/2$ copies de $GL(d_{\rho},F)$.

Il faut maintenant comparer les facteurs de transfert: on a $\gamma_{G}=\times_{i\in [1,(a-a_{-})/2} m_{i}\times \gamma'$ et $\gamma_{H}=_{i\in [1,(a-a_{-})/2} m_{i}\times \gamma'_{1}\times \gamma_{2}$ avec chaque $m_{i}$ dans un groupe $GL(d_{\rho},F)$, $\gamma'\in SO(2n-d_{\rho}(a-a_{-})+1,F)$, $\gamma'_{1}\in SO(2n_{1}-d_{\rho}(a-a_{-})+1,F)$ et $\gamma_{2}\in SO(2n_{2}+1,F)$. Et on est dans la situation o\`u la classe stable de $\gamma'_{1}\times \gamma_{2}$ correspond \`a la classe stable de 
$\gamma'$ dans l'endoscopie pour $SO(2n-d_{\rho}(a-a_{-})+1,F)$. On v\'erifie que $$\Delta(\gamma_{G},\gamma_{H})=\Delta(\gamma',\gamma'_{1}\gamma_{2})$$
 le premier facteur de transfert \'etant pour $SO(2n+1,F)$ et le second pour $SO(2n-d_{\rho}(a-a_{-})+1,F)$; cela r\'esulte du calcul des int\'egrales orbitales pour les \'el\'ements d'un Levi.

Finalement on \'etablit que la distribution sur $M_{a,a_{-}}$
$$
\sum_{\pi\in \Pi(\psi)}\epsilon_{a}(\pi)\epsilon_{a_{-}}(\pi)c_{\pi}proj_{\rho\vert\,\vert^{(a_{-}+1)/2}}\cdots proj_{\rho\vert\,\vert^{(a-1)/2}}res_{M_{a,a_{-}}}\pi$$
est un transfert de la distribution stable sur $M^1_{a,a_{-}}$:
$$
\biggl(\sum_{\pi_{1}\in \Pi(\psi_{1})}proj_{\rho\vert\,\vert^{(a_{-}+1)/2}}\cdots proj_{\rho\vert\,\vert^{(a-1)/2}}res_{M^1_{a,a_{-}}}\pi_{1}\biggr)\biggl(\sum_{\pi_{2}\in \Pi(\psi_{2})}c_{\pi_{2}}tr\, \pi_{2}\biggr).
$$

Le r\'esultat principal de \cite{ams} d\'ecrit exactement les termes $proj_{\rho\vert\,\vert^{(a_{-}+1)/2}}\cdots proj_{\rho\vert\,\vert^{(a-1)/2}}res_{M_{a,a_{-}}}\pi$  quand $\pi$ parcourt $\Pi(\psi)$. Exactement cette trace est $0$ pour la moiti\'e des repr\'esentations pr\'ecis\'ement celles pour lesquelles notre caract\`ere $\epsilon_{\pi}$ associ\'e \`a tout \'el\'ement de $\Pi(\psi)$\cite{europe}, \cite{ams} par  v\'erifie $\epsilon_{\pi}(\rho,a)\neq \epsilon_{\pi}(\rho,a_{-})$ 
et pour les autres est une repr\'esentation irr\'eductible de la forme $$\rho\vert\,\vert^{(a-1)/2}\otimes \cdots \otimes \rho\vert\,\vert^{(a_{-}+1)/2}\otimes \pi'$$ o\`u $\pi'$ est l'un des 2 sous-modules 
irr\'eductibles de l'induite $St(\rho,a_{-})\times \pi''$, avec $\pi''$ est un \'el\'ement de $\Pi(\psi_{2})$; si $a_{-}=0$, il n'y a pas d'induite et on obtient directement $\pi''$. 
On note simplement $Jac_{(a-1)/2, \cdots, (a_{-}+1)/2}\pi$ cette repr\'esentation irr\'eductible. On a 
d\'ecrit $\Pi(\psi_{1})$ et avec cela on v\'erifie que:
$$\biggl(\sum_{\pi_{1}\in \Pi(\psi_{1})}proj_{\rho\vert\,\vert^{(a_{-}+1)/2}}\cdots proj_{\rho\vert\,\vert^{(a-1)/2}}res_{M^1_{a,a_{-}}}\pi_{1}\biggr)$$
est $\otimes_{j\in [(a-1)/2,(a_{-}+1)/2]}\rho\vert\,\vert^{j}\otimes ind St(\rho,a_{-})$ o\`u l'induite est pour le groupe $SO(2a_{-}+1,F)$ \`a partir du parabolique de Levi $GL(a_{-},F)$. Ainsi la distribution pour $SO(2n-d_{\rho}(a-a_{-})+1,F)$ 
$$
\sum_{\pi\in \Pi(\psi)}\epsilon_{a}(\pi)\epsilon_{a_{-}}(\pi)c_{\pi}Jac_{(a-1)/2, \cdots, (a_{-}+1)/2}\pi\eqno(5)
$$
est un transfert de la distribution (rappelons que $n_{1}=(a+a_{-})/2$):
$$
\biggl(ind_{GL(a_{-}d_{\rho},F)}^{SO(2a_{-}+1,F)}\, St(\rho,a_{-})\biggr)\otimes \biggl(\sum_{\pi_{2}\in \Pi(\psi_{2})}c_{\pi_{2}}\pi_{2}\biggr).\eqno(6) 
$$

Si $a_{-}=0$, on peut tout de suite conclure:  les groupes pour (5) et (6) sont les m\^emes et le transfert est l'identit\'e. Donc en particulier les coefficients $c_{\pi_{2}}$ et $c_{\pi}$ \'etant positifs on doit aussi avoir $\epsilon_{a}(\pi)=1$ si $Jac_{(a-1)/2, \cdots, 1/2}\pi\neq 0$. De plus comme on l'a rappel\'e, la moiti\'e des \'el\'ements de $\Pi(\psi)$ ont cette propri\'et\'e de non nullit\'e. On calcule le cardinal du groupe des caract\`eres de $Cent_{Sp(2n,{\mathbb C})}(\psi)/Cent(Sp(2n,{\mathbb C}))$ triviaux sur $Z_{\rho\otimes \sigma_a}$; ce groupe est facilement mis en bijection avec le groupe des caract\`eres de  
$Cent_{Sp(2n-ad_{\rho},{\mathbb C})}(\psi_{2})/Cent(Sp(2n-ad_{\rho},{\mathbb C}))$. 
Le cardinal de ce groupe est donc $2^{\ell_{\psi_{2}}-1}$ ce qui est la moiti\'e du cardinal de $\Pi(\psi)$. Comme on a d\'ej\`a montr\'e que l'application $\pi\in \Pi(\psi)$ associe $\epsilon_{\cal A}$ est une bijection sur l'ensemble des caract\`eres du centralisateur de $\psi$ triviaux sur le centre de $Sp(2n,{\mathbb C})$ (cf. \ref{classification}), on voit que la propri\'et\'e $\epsilon_{a}(\pi)=1$ caract\'erise les repr\'esentations $\pi$ de $\Pi(\psi)$ telles que $Jac_{(a-1)/2, \cdots, 1/2}\pi\neq 0$ ce qui est l'\'enonc\'e cherch\'e dans ce cas.

\

Supposons maintenant que $a_{-}>0$. Soient $\pi,\pi'\in \Pi(\psi)$ et $\pi_{2}\in \Pi(\psi_{2})$ tels que l'on ait l'\'egalit\'e:
$$
Jac_{(a-1)/2, \cdots, (a_{-}+1)/2}\pi\oplus Jac_{(a-1)/2, \cdots, (a_{-}+1)/2}\pi'=St(\rho,a_{-})\times \pi_{2}.
$$
On montre d'abord que $c_{\pi}=c_{\pi'}=c_{\pi_{2}}$; c'est un probl\`eme de transfert stable. D'une part on sait que
$\sum_{\pi\in \Pi(\psi)}c_{\pi}\pi$ est le transfert d'un prolongement \`a $\tilde{GL}(2n,F)$ de la repr\'esentation temp\'er\'ee $\pi(\psi)$ de $GL(2n,F)$ associ\'ee \`a $\psi$. D'autre part on sait que $\sum_{\pi_{2}\in \Pi(\psi_{2})}c_{\pi_{2}}\pi_{2}$ est le transfert stable d'un prolongement \`a $\tilde{GL}(2n_{2},F)$ de la repr\'esentation temp\'er\'ee $\pi(\psi_{2})$ de $GL(2n_{2},F)$ associ\'ee \`a $\psi_{2}$. On sait que $\pi(\psi)$ est l'unique sous-module irr\'eductible de l'induite
$$
\theta\vert\,\vert^{(a-1)/2}\times \cdots \times \rho\vert\,\vert^{(a_{-}+1)/2}\times St(\rho,a_{-}) \times \pi(\psi_{2}) \times St(\rho,a_{-}) \times \rho\vert\,\vert^{-(a_{-}+1)/2}\times \cdots \times \rho\vert\,\vert^{-(a-1)/2}.
$$
Donc on passe de $\pi(\psi)$ \`a $St(\rho,a_{-})\times \pi(\psi_{2}) \times St(\rho,a_{-})$ en prenant des modules de Jacquet (cf. \cite{transfert} 4.2) et on passe de $\pi(\psi_{2})$ \`a la m\^eme repr\'esentation en induisant. Cela donne imm\'ediatement l'\'egalit\'e
$$
\sum_{\pi\in \Pi(\psi)}c_{\pi}Jac_{(a-1)/2, \cdots, (a_{-}-1)/2}\pi=\sum_{\pi_{2}\in \Pi(\psi_{2})}c_{\pi_{2}}(St(\rho,a_{-})\times \pi_{2}).
$$
On en d\'eduit l'assertion en d\'ecomposant chaque repr\'esentation comme d\'ej\`a expliqu\'e. Pour $\pi_{2}\in \Pi(\psi_{2})$, on note $\Pi(\pi_{2})$ les 2 repr\'esentations, $\pi'$, $\pi''$  de $\Pi(\psi)$ telles que $$Jac_{(a-1)/2, \cdots, (a_{-}+1)/2}(\pi'\oplus \pi'')=St(\rho,a_{-})\times \pi_{2}.$$

On a donc sur $SO(2n_{2}+2d_{\rho}a_{-}+1,F)$ la distribution:
$$
\sum_{\pi_{2}\in \Pi(\psi_{2})}c_{\pi_{2}}\sum_{\pi\in \Pi(\pi_{2})}\epsilon_{a}(\pi)\epsilon_{a_{-}}(\pi)Jac_{(a-1)/2, \cdots, (a_{-}+1)/2} \pi
$$
qui est un transfert de la distribution sur $SO(2d_{\rho}a_{-}+1,F)\times SO(2n_{2}+1,F)$:
$$
\biggl(ind_{GL(d_{\rho}a_{-},F)}^{SO(2d_{\rho}a_{-}+1,F)}St(\rho,a_{-})\biggr)\biggl( \sum_{\pi_{2}}c_{\pi_{2}}\pi_{2}\biggr).
$$
On \'ecrit encore l'\'egalit\'e qui en r\'esute pour tout \'el\'ement $\gamma\in GL(d_{\rho}a_{-})\times SO(2n_{2}+1,F)$ dont la partie dans $GL(d_{\rho}a_{-})$ est elliptique. Ensuite on projette sur le sous-espace vectoriel de $I_{cusp}(GL(d_{\rho}a_{-}))$ image de $St(\rho,a_{-})$. On v\'erifie ais\'ement que les \'el\'ements de $\Pi(\psi_{2})$ n'ont pas de module de Jacquet qui peuvent contribuer \`a une telle projection. Le facteur de transfert  $\Delta(\gamma,\gamma)$, o\`u $\gamma$ est d'abord vu comme un \'el\'ement de $SO(2n_{2}+2d_{\rho}a_{-}+1,F)$ puis du groupe endoscopique $SO(2d_{\rho}a_{-}+1,F)\times SO(2n_{2}+1,F)$ vaut 1.

On obtient alors une \'egalit\'e:
$$
\sum_{\pi_{2}\in \Pi(\psi_{2})}c_{\pi_{2}}\biggl(\sum_{\pi\in \Pi(\pi_{2})}(\epsilon_{a}(\pi)\epsilon_{a_{-}}(\pi)\biggr)\pi_{2}=
$$
$$
2\times \biggl(\sum_{\pi_{2}\in \Pi(\psi_{2})} c_{\pi_{2})}\pi_{2}\biggr).
$$
D'o\`u $\epsilon_{a}\epsilon_{a_{-}}(\pi)=1$ pour tout $\pi_{2}$ et tout $\pi\in \Pi(\pi_{2})$. Ensuite on conclut comme dans le cas o\`u $a_{-}=0$; la moiti\'e des \'el\'ements, $\pi$,  de $\Pi(\psi)$ v\'erifie \`a la fois $Jac_{(a-1)/2, \cdots, (a_{-}+1)/2}\pi \neq 0$ et $\epsilon_{a}(\pi)\epsilon_{a_{-}}(\pi)=1$. Mais chacune des 2 conditions est satisfaite par exactement la moiti\'e des \'el\'ements de $\Pi(\psi)$, elles sont donc \'equivalentes. C'est l'assertion cherch\'ee.
\subsection{Description des param\`etres des repr\'esentations cuspidales\label{cuspidal}}
\bf Th\'eor\`eme. \sl L'ensemble des repr\'esentations cuspidales de $SO(2n+1,F)$ est en bijection avec l'ensemble des couples $(\psi,\epsilon)$ o\`u $\psi$ est un morphisme sans trou de $W_{F}\times SL(2,{\mathbb C})$ dans $Sp(2n,{\mathbb C})$ et $\epsilon$ est un caract\`ere altern\'e du centralisateur de $\psi$ de restriction trivial au centre de $Sp(2n,{\mathbb C})$ (on renvoit \`a l'introduction pour une pr\'ecision sur les notations).
\rm

\

Soit $\pi$ une repr\'esentation cuspidal et $\psi$ le morphisme tel que $\pi\in \Pi(\psi)$; on a d\'ej\`a v\'erifi\'e que $\psi$ est sans trou; on note $\epsilon_{\cal A}(\pi)$ le caract\`ere associ\'e par Arthur \`a $\pi$; il est altern\'e par \ref{caractereetmoduledejacquet}. R\'eciproquement soit $\psi$ un morphisme sans trou et $\pi\in \Pi(\psi)$. Soit $\rho$ une repr\'esentation cuspidale irr\'eductible 
et unitaire d'un $GL(d_{\rho},F)$, $x\in {\mathbb R}$ et $\sigma$ une repr\'esentation irr\'eductible de $SO(2n-d_{\rho}+1,F)$ tel que $$\pi \hookrightarrow \rho\vert\,\vert^x\times \sigma.\eqno(1)$$Le fait que $\pi$ soit cuspidal est exactement \'equivalent \`a ce qu'il ne soit pas possible de trouver de telles donn\'ees satisfaisant (1). Supposons que $\pi$ n'est pas cuspidal et montrons que $\epsilon_{\cal A}\pi$ n'est pas altern\'e. Cela suffira. On fixe donc $\rho,x$ satisfaisant (1); $\sigma$ ne joue pas de 
r\^ole. On sait a priori que $x$ est un demi-entier (cf. \ref{integralite}) et d'apr\`es \cite{europe} par exemple qu'en notant encore $\rho$ la repr\'esentation irr\'eductible de $W_{F}$ correspondant \`a $\rho$ la repr\'esentation d\'efinie par $\psi$ de $W_{F}\times SL(2,{\mathbb C})$ contient comme sous-repr\'esentation la repr\'esentation $\rho\otimes \sigma_{2x+1}$; de plus comme $\pi$ est une 
s\'erie discr\`ete, on a s\^urement $x>0$. On \'ecrit $a=2x+1$ et $a\geq 2$. Comme $\psi$ est sans trou, $a_{-}$ est d\'efini et vaut $a-2$. Il suffit d'appliquer \ref{caractereetmoduledejacquet} pour voir que $\epsilon_{\cal A}(\pi)$ n'est pas altern\'e.

 \subsection{Au sujet de la combinaison lin\'eaire stable dans un paquet.\label{calculdescoefficients}}
Soit encore $\psi$ un morphisme $\theta$-stable et $\theta$ discret de $W_{F}\times SL(2,{\mathbb C})$ dans $Sp(2n,{\mathbb C})$. On note $\sum_{\pi\in \Pi(\psi)}c_{\pi}\pi$ la distribution stable engendr\'ee par les \'el\'ements du paquet $\Pi(\psi)$ avec $c_{\pi}$ des r\'eels positifs; pour fixer vraiment les $c_{\pi}$ on demande \`a cette distribution d'\^etre le transfert de la trace de la repr\'esentation temp\'er\'ee $\pi(\psi)$ de $GL(2n,F)$ prolong\'e \`a $\tilde{GL}(2n,F)$. On remarque que $\psi$ est uniquement 
d\'etermin\'e par $\pi$ et que $c_{\pi}$ est donc aussi uniquement d\'etermin\'e par $\pi$.

 Il est naturel de conjecturer que les $c_{\pi}$ sont tous \'egaux \`a 1. Ceci est d\'emontr\'e dans \cite{waldspurger2} si la restriction de $\psi$ \`a $W_{F}$ est triviale sur le groupe de ramification mod\'er\'ee. Par changement de base, on devrait pouvoir \'etendre le r\'esultat de Waldspurger \`a tous les morphismes tels que la restriction de $\psi$ \`a $W_{F}$ se factorise par le groupe de Weil 
 d'une extension r\'esoluble de $F$. Mais dans l'\'etat actuel, il semble difficile d'aller au del\`a et la remarque suivante est un substitut qui permettra d'\'ecrire les distributions stables (et endoscopiques) pour tous les paquets pas seulement ceux qui sont temp\'er\'es. Pour cela, il faut rappeler ce 
 qu'est le support cuspidal partiel d'une repr\'esentation irr\'eductible, $\pi$ d'un groupe classique, ici $SO(2n+1,F)$: c'est l'unique repr\'esentation cuspidale, $\pi_{cusp}$, d'un groupe de la forme $SO(2n_{cusp}+1,F)$ (ce qui d\'efinit l'entier $n_{cusp}$) tel que $\pi$ soit un sous-quotient irr\'eductible d'une induite de la forme $\sigma\times \pi_{cusp}$ o\`u $\sigma$ est une repr\'esentation 
 irr\'eductible de $GL(n-n_{cusp},F)$. Comme une repr\'esentation cuspidale, comme $\pi_{cusp}$, est dans un unique paquet de s\'erie discr\`ete, on a bien d\'efini $c_{\pi_{cusp}}$ tout comme on a d\'efini $c_{\pi}$ ci-dessus.

\

\bf Remarque. \sl Soit $\pi$ une s\'erie discr\`ete dont on note $\pi_{cusp}$ le support cuspidal partiel, alors $c_{\pi}=c_{\pi_{cusp}}$.\rm

\

On le d\'emontre par r\'ecurrence sur $n$ et on r\'eutilise la d\'emonstration de \ref{caractereetmoduledejacquet}. Dans la preuve de \ref{caractereetmoduledejacquet} on a d\'ej\`a calcul\'e $c_{\pi}$ en fonction d'une s\'erie 
discr\`ete d'un groupe $SO(2n'+1,F)$ avec $n'<n$ sous l'hypoth\`ese que la repr\'esentation 
d\'efinie par $\psi$ contient une sous-repr\'esentation $\rho\otimes \sigma_a$ avec $a_{-}$ d\'efini et $\pi$ v\'erifiant $\epsilon_{a}(\pi)=\epsilon_{a_{-}}(\pi)$ avec les notations de cette preuve. Il est imm\'ediat de voir que $\pi_{cusp}$ est aussi le support cuspidal partiel de la repr\'esentation $\pi_{2}$ 
de $SO(2n_{2}+1,F)$ pour laquelle on a montr\'e que $c_{\pi}=c_{\pi_{2}}$. D'o\`u le r\'esultat par 
r\'ecurrence dans ce cas. On est donc ainsi ramen\'e au cas o\`u $\epsilon_{\cal A}(\pi)$ (le caract\`ere associ\'e par \cite{arthurnouveau} 30.1 \`a $\pi$) est altern\'e. Faisons cette hypoth\`ese sur $\epsilon_{\cal A}(\pi)$. 
D'apr\`es \ref{classification} et la description des 
param\'etres des repr\'esentations cuspidales (\ref{cuspidal}), $\pi$ est cuspidal si et seulement si $\psi$ est sans trou; il n'y a donc rien \`a d\'emontrer dans le cas o\`u $\psi$ est sans trou. 

Il faut donc voir encore le cas o\`u $\psi$ a des trous, ou encore, le cas o\`u 
il existe $\rho\otimes \sigma_a$ une sous-repr\'esentation irr\'eductible de $\psi$ avec $a>2$ telle que $\rho\otimes \sigma_{a-2}$ ne soit pas une sous-repr\'esentation de $\psi$. On fixe une telle 
repr\'esentation $\rho\otimes \sigma_a$ et on note $d_{\rho}$ la dimension de la repr\'esentation $\rho$ et on notera aussi $\rho$ la repr\'esentation cuspidale de $GL(d_{\rho},F)$ associ\'ee \`a $\rho$ 
par la correspondance locale de Langlands (\cite{harris},\cite{henniart}). On note $\psi'$ le morphisme de $W_{F}\times SL(2,{\mathbb C})$ dans $Sp(2(n-d_{\rho}),{\mathbb C})$ qui 
est la somme des sous-repr\'esentations irr\'eductibles incluses dans $\psi$ sauf $\rho\otimes \sigma_a$ qui est remplac\'ee par $\rho\otimes \sigma_{a-2}$. On a montr\'e en \cite{europe} et \cite{ams} 
(cf. l'introduction de \cite{ams}) que l'application $Jac_{(a-1)/2}$ \'etablit une bijection entre $\Pi(\psi)$ et $\Pi(\psi')$; rappelons ce qu'est cette application. Soit $\pi\in \Pi(\psi)$, alors il existe une 
unique 
repr\'esentation irr\'eductible $\pi'$ de $SO(2(n-d_{\rho})+1,F)$ telle que $\pi$ soit un sous-module 
de l'induite $\rho\vert\,\vert^{(a-1)/2}\times \pi'$. Et on a montr\'e que $\pi'\in \Pi(\psi')$ et que 
tous les 
\'el\'ements de $\Pi(\psi')$ sont obtenus ainsi exactement une fois.

Par la compatibilit\'e du transfert stable \`a la prise de module de Jacquet (\cite{transfert} 4.2), on obtient le fait que 
$
\sum_{\pi\in \Pi(\psi)}c_{\pi} Jac_{(a-1)/2}\pi$ est un transfert d'un prolongement (bien d\'etermin\'e) de la repr\'esentation temp\'er\'ee $\pi(\psi')$ de $GL(2(n-d_{\rho}),F)$ \`a $\tilde{GL}(2(n-d_{\rho}),F)$. D'o\`u l'\'egalit\'e $c_{\pi}=c_{Jac_{(a-1)/2}\pi}$. Et on applique l'hypoth\`ese de r\'ecurrence \`a $Jac_{(a-1)/2}\pi$. On obtient l'\'egalit\'e $c_{\pi}=c_{\pi_{cusp}}$ puisqu'il est bien clair que le support cuspidal partiel de $\pi$ est le m\^eme que celui de $Jac_{(a-1)/2}\pi$. Cela termine la preuve.

\section{Le cas des groupes orthogonaux impairs non d\'eploy\'es\label{casnondeploye}}
Ici $G(n)$ est un groupe orthogonal d'un espace de dimension $2n+1$ d'une forme orthogonale ayant un noyau anisotrope de dimension 3. Les preuves ci-dessous n'ont rien d'originales et s'appliquent beaucoup plus g\'en\'eralement pour passer d'un groupe quasid\'eploy\'e \`a une forme int\'erieure.

Et on note $G_{d}(n)$ la forme 
d\'eploy\'ee de $SO(2n+1,F)$. On peut encore d\'efinir comme en \ref{argumentlocal}, l'espace $I^{st}_{cusp}(G(n))$. Soit $\pi$ une s\'erie discr\`ete, l'argument donn\'e en loc.cit. s'applique pour montrer que l'image du caract\`ere de $\pi$ dans $I^{st}_{cusp}(G(n))$ est non nul. Il existe donc une distribution stable de la forme
$$
\sum_{\pi'\in {\cal P}}c_{\pi'}tr\, \pi',\eqno(1)
$$
o\`u ${\cal P}$ est un ensemble de repr\'esentations elliptiques de $G(n)$ et $c_{\pi'}$ sont des \'el\'ements de ${\mathbb C}$ et tel que $\pi\in {\cal P}$.

D'apr\`es \cite{selecta} 3.5 (un peu amelior\'e en \cite{inventiones} 4.5),  le transfert induit un isomorphisme de $I^{st}_{cusp}(G(n))$ sur $I^{st}_{cusp}(G_{d}(n))$. Ainsi, il existe un ensemble ${\cal P}_{d}$ de repr\'esentations elliptiques de $G_{d}(n)$ et des nombres complexes, $c_{\pi'_{d}}$ pour $\pi'_{d}\in {\cal P}_{d}$ tels que $$\sum_{\pi'_{d}\in {\cal P}_{d}}c_{\pi'_{d}}tr\pi'_{d}\eqno(2)$$
soit un transfert stable de (1). R\'eciproquement, \'etant donn\'e un paquet stable pour $G_{d}(n)$, c'est-\`a-dire une combinaison lin\'eaire \`a coefficients complexes de caract\`eres, stable, il existe une combinaison lin\'eaire \`a coefficient complexes de caract\`eres pour $G(n)$ qui est stable et qui se transf\`ere en la combinaison de d\'epart. Les caract\`eres sont les caract\`eres des 
repr\'esentations a priori elliptiques.

\

Ainsi pour $\pi$ une s\'erie discr\`ete de $G(n)$, il existe au moins un morphisme $\psi$ de $W_{F}\times SL(2,{\mathbb C})$ dans $G^*(n)$ tel que $\pi$ soit dans un paquet stable se tranf\'erant en le paquet $\Pi(\psi)$ form\'e de s\'eries discr\`etes de $G_{d}(n)$. Ceci suppose \'evidemment les hypoth\`eses 
d\'ej\`a faites au sujet des lemmes fondamentaux. On dira encore que $\pi\in \Pi(\psi)$, $\Pi(\psi)$ \'etant maintenant vu comme un ensemble de repr\'esentations elliptiques de $G(n)$.

\ 

\bf Th\'eor\`eme 1. \sl  Soit $\pi_{0}$ une repr\'esentation cuspidale de $G(n_{0})$ et $\rho$ une repr\'esentation cuspidale autoduale irr\'eductible d'un $GL(d_{\rho},F)$. Soit $x_{0}\in {\mathbb R}$ tel que l'induite $\rho\vert\,\vert^{x_{0}}\times \pi_{0}$ soit r\'eductible. Alors $x_{0}\in 1/2 {\mathbb Z}$.\rm

\

On fixe $\pi_{0}$ et $x_{0}$ comme dans l'\'enonc\'e; si $x_{0}=0$, il n'y a rien \`a d\'emontrer et on suppose donc que $x_{0}>0$. On note ici $\pi$ l'unique sous-module irr\'eductible de $\rho\vert\,\vert^{x_{0}}\times \pi$. On reprend les notations, $Jac_{x}$ de \ref{integralite}. On applique $Jac_{x_{0}}$ \`a (1);  soit  $\pi'\in {\cal P}$ tel que $Jac_{x_{0}}\pi'\neq 0$ et contient $\pi_{0}$ dans sa d\'ecomposition en irr\'eductible. Cela entra\^{\i}ne que $\pi'$ est un sous-quotient irr\'eductible de l'induite $\rho\vert\,\vert^{x_{0}}\times \pi_{0}$ et comme $\pi'$ est elliptique cela n\'ecessite que $\pi'=\pi$. Ainsi $Jac_{x_{0}}$ appliqu\'e \`a (1) est de la forme:
$$
c_{\pi}\pi_{0}\oplus \tau \eqno(3)
$$
o\`u $\tau$ dans (3) est une combinaison lin\'eaire de repr\'esentations irr\'eductibles dont aucun n'est \'equivalente \`a $\pi_{0}$. On v\'erifie que $Jac_{x_{0}}$ appliqu\'e \`a (2) est un transfert de (3). En particulier, $Jac_{x_{0}}$ appliqu\'e \`a (2) n'est pas nul. Ou encore, il existe une repr\'esentation elliptique de $G_{d}(n)$, $\pi_{d}$, tel que $Jac_{x_{0}}\pi_{d}\neq 0$. Mais comme on sait que le support cuspidal des repr\'esentations elliptiques de $G_{d}(n)$ est comme celui des s\'eries discr\`etes, demi-entier, on en 
d\'eduit que $x_{0}$ est demi-entier. D'o\`u le th\'eor\`eme.

\

On d\'emontre de la m\^eme fa\c{c}on les th\'eor\`emes ci-dessous:

\

\bf Th\'eor\`eme 2. \sl Avec les hypoth\`eses et notations du th\'eor\`eme pr\'ec\'edent en particulier $\pi_{0}$ est une repr\'esenta\-tion cuspidale de $G(n)$. Soit $\psi_{0}$ un morphisme de $W_{F}\times SL(2,{\mathbb C})$ tel que $\pi_{0}\in \Pi(\psi_{0})$. Alors $\psi_{0}$ est $\theta$-discret, sans trou et les points de r\'eductibilit\'e pour les induites de la forme $\rho\vert\,\vert^x\times \pi_{0}$ se calculent comme dans le cas d\'eploy\'e, c'est-\`a-dire que l'on a avec les notations de ce cas:
$$
x_{\rho,\pi_{0}}=(a_{\rho,\psi_{0}}+1)/2.
$$

\bf Th\'eor\`eme 3. \sl Soit $\pi$ une s\'erie discr\`ete irr\'eductible de $G(n)$ et soit $\psi$ un morphisme de $W_{F}\times SL(2,{\mathbb C})$ dans $Sp(2n,{\mathbb C})$ tel que $\pi\in \Pi(\psi)$. Alors $\psi$ est $\theta$-discret. De plus ce morphisme $\psi$ est uniquement d\'etermin\'e par les blocs de Jordan de $\pi$ comme dans le cas d\'eploy\'e et correspond \`a celui qui a \'et\'e associ\'e \`a $\pi$ par \cite{europe} et \cite{ams}

\
 \rm
Si $\psi$ n'est pas $\theta$-discret, comme dans le cas d\'eploy\'e toute repr\'esentation de $\Pi(\psi)$ est sous-module d'une induite convenable. Ceci est exclu pour $\pi$ qui est une s\'erie discr\`ete, d'o\`u le th\'eor\`eme 3. La fin du th\'eor\`eme se d\'emontre comme dans le cas d\'eploy\'e, nos r\'ef\'erences, \cite{algebra}, \cite{europe} et \cite {ams} ne font pas l'hypoth\`ese que le groupe est d\'eploy\'e.

\

Pour terminer remarquons que l'on n'a pas ici calcul\'e le nombre d'\'el\'ements d'un paquet $\Pi(\psi)$ et par voie de cons\'equence obtenu une param\'etrisation des repr\'esentations cuspidales de $G(n)$; la m\'ethode du cas d\'eploy\'e s'appuie sur les r\'esultats de transfert endoscopique d'Arthur qu'il faudrait \'etendre. Il y a sans doute une m\'ethode locale comme pour le transfert stable et cela vaudrait la peine de l'\'ecrire en g\'en\'eral.

\end{document}